\documentclass[a4paper,fleqn]{cas-sc}

\usepackage[authoryear,longnamesfirst]{natbib}

\usepackage{todonotes}
\usepackage{amsmath}
\usepackage{siunitx}
\usepackage[export]{adjustbox}

\usepackage[final]{listings}
\usepackage{tcolorbox}
\tcbuselibrary{skins,breakable}
\newtcolorbox{codebox}[2][]{
    colback=white,
    colframe=blue,
    coltitle=black,
    breakable = true,
    boxsep = -2.5mm,
    fonttitle=\small\ttfamily,enhanced,
    attach boxed title to top right={xshift = -0.5cm, yshift = -5.5mm},
    boxed title style={colback=white, colframe=blue},
    title=#2,#1
}

\newcommand{\dealii}{{\texttt{deal.II}}}

\begin{document}
\let\WriteBookmarks\relax
\def\floatpagepagefraction{1}
\def\textpagefraction{.001}
\shorttitle{Parallel flow routing}
\shortauthors{W. Bangerth}

\title [mode = title]{Massively parallel flow routing and drainage area determination}

\author{Wolfgang Bangerth}[orcid=0000-0003-2311-9402]
\ead{bangerth@colostate.edu}

\credit{Conceptualization of this study, Methodology, Software, Writing}

\affiliation{organization={Department of Mathematics, Department of
    Geosciences, Colorado State University},
                addressline={1874 Campus Delivery}, 
                city={Fort Collins},
                postcode={80524}, 
                state={CO},
                country={United States}}

\begin{abstract}
Digital elevation models (DEMs) have reached resolutions and sizes that only parallel computaters can efficiently process. One important application of DEMs is predicting how much water flows where, the so-called ``flow routing problem'' (a variation of which is the problem of determining the drainage area upstream of a point in a DEM). The traditional algorithm for flow routing is sequential, and attempts to parallelize this method have so far only been moderately successful. Herein, we build on earlier work in \cite{Richardson2014} and propose an algorithm and several variations that can efficiently solve the flow routing problem on very large models with very large numbers of parallel processes. For the largest model we use, with 1.88 billion points, the best algorithm herein can route water in 4.0 seconds on \num{12288} processes of a computer cluster.
\end{abstract}


\maketitle

\section{Introduction}
\label{sec:introduction}

Flow routing refers to the problem of determining how much water flows at every point of a landscape, given a digital elevation model (DEM) and knowledge of how much water originates where, for example through rainfall. (The related problem of determining drainage area can be stated as seeking the amount of water that results from a spatially constant rainfall rate because in that case, the amount of water is proportional to the drainage area.) 
Flow routing has traditionally been solved using algorithms that either work their way upward in drainages, or work their way downward from mountain tops. The intellectual basis for the work herein is an algorithm that is a model of elegance because it is not non-obvious yet simple to describe: In a first phase, one sorts all grid points in the DEM from high to low; in a second phase, points are processed in this high-to-low order in two steps: (i) The amount of water available at a point equals the amount of rainfall in the area that corresponds to this point, plus whatever water it has received from previously processed points; (ii) the available water is then, depending on modeling assumptions, given either to the lowest-lying neighbor of the current point, or partitioned among its lower-lying neighbors. This algorithm works because the order of processing ensures that when a point in the DEM is visited, all higher-lying neighbors have already been processed, and because a point can only receive water from higher-lying neighbors, we have all the knowledge we need about how much water is available at the current point. There can be no cycles, and we never have to visit points a second time, or do other kinds of iterations.

For a digital elevation model with $N$ points, the sorting phase requires ${\cal O}(N\, \log_2 N)$ operations, whereas the processing phase touches each point exactly once and so requires ${\cal O}(N)$ operations. It is difficult to imagine a method that would beat this algorithm, either in speed or elegance. Yet, the algorithm is unusable in today's environment where digital elevation models can have $N=10^{11}$ points or more \cite{ASTER}:
For models of this size, only algorithms running on many machines in parallel can be meaningfully used to process these amounts of data; on the other hand, the one outlined in the previous paragraph is inherently sequential because it visits points in a specific order one after the other.

In order to allow for flow routing in large elevation models, there have been attempts to parallelize the algorithm above on shared-memory multi-core machines, by keeping track of which points have to be processed before which others. This can be done through parallel priority queues, or by counting how many unmet dependencies each point still has, but these ideas do not lead to speed-ups that really make a difference for the largest data sets. Similarly, attempts have been made to utilize distributed memory machines and the Message Passing Interface (MPI), but these also were not very encouraging: To give just one example, \cite{wallace2010parallel} report a speed-up of around 3 when using 8 processes, and no more than 4 when using up to 128 processes.%
\footnote{Herein, I use the term ``process'' to describe a single program running on a machine, as is common in the literature on parallel programming using the Message Passing Interface (MPI). A parallel program then consists of multiple processes, each of which has its own memory space. Oftentimes, one runs one process per processor of a cluster or (on today's multi-core machines) one process per processor core.} Section 2.4 of \cite{Richardson2014} provides a nice overview of many of the attempts made, along with the unsatisfactory state of affairs.

Moreover, current digital elevation models exceed the memory size of even large workstations, making them unsuitable for algorithms that require access to all data at once (i.e., ``shared memory'' algorithms):  algorithms that require having the entire data set in memory on a single machine can no longer be used at the scale of today's DEMs. 

As a consequence of these considerations, algorithms that work with shared memory, or algorithms that only scale to a handful or perhaps a few dozen processors will not be sufficient in the future. There clearly is a need to look at the problem in an entirely different way.

One basis for such an alternative perspective is the realization -- apparently independently proposed in \cite{Eddins2007}, \cite{Schwanghart2010}, and Richardson, Hill, and Perron (see \cite[Section 3.1]{Richardson2014} for the claim of independent discovery) -- that flow routing and drainage area determination can be formulated as solving a system of linear equations. The main contribution of \cite{Richardson2014} is that there are efficient algorithms for solving systems of linear equations that can be parallelized and that can be run on large clusters of machines that do not share access to common memory. The algorithm proposed therein is called Implicit Drainage Area (IDA); the paper does not compare run-times between the best algorithm on a single process and their parallel implementation, but the authors demonstrate a speed-up of about 3 when going from 48 to 192 processes. That said, reading between the lines and taking into account (i) the fact that their algorithm requires hundreds of iterations (see their Fig.~13) and (ii) the fact that each iteration costs about as much as a single high-to-low solve, one can conjecture that the IDA algorithm is perhaps slightly faster, but likely not substantially so, on 48 cores than the high-to-low method on a single core.

The contribution of the current paper is to show that the formulation of flow routing and drainage area determination as a linear system can be used as the basis for a much better algorithm than IDA using insight into which algorithms and preconditioners are suitable for the specific type of linear system we have here. In particular, in the remainder of this paper I will present an algorithm and several variations with the following properties:
\begin{enumerate}
    \item If run on a single process, the methods are all equivalent to the traditional high-to-low sweep, and consequently just as fast.
    \item The algorithms can be run on large numbers of processes, showing a speed-up of 52 compared to the single process case when run on 128 processes using a digital elevation model with 177 million points, and solving a problem with 1.88 billion points in 4.0 seconds on \num{12288} processes that could not at all be solved with fewer than 8 processes.
    \item The algorithms partition memory requirements equally among all processes, with every process only requiring knowledge of that part of the the digital elevation model that corresponds to those points it ``owns'' plus one layer of neighboring points. As a consequence, the methods can efficiently deal with very large DEM data sets that do not fit into the memory of a single machine, but that can be loaded into the memory spaces of hundreds or thousands of machines taken together.
\end{enumerate}

The remainder of this paper is structured as follows: In Section~\ref{sec:basic-algorithm}, I will describe the classical ``high-to-low'' algorithm and then outline the basic idea of the parallel flow routing algorithm. Section~\ref{sec:optimized-algorithm} then covers optimizations that lead to identical results but reduce memory usage or run time, and simplify the set-up of the algorithm; the section then also reformulates the problem further, allowing for a more sophisticated preconditioner. I will then show experimental results with this algorithm on some large DEM data sets that have up to 1.88 billion points and using thousands of parallel processes in Section~\ref{sec:evaluation}; these results will show that one can attain substantial speed-ups from parallel computing when using sufficiently large models. I will conclude and consider where this line of research might lead in Section~\ref{sec:conclusions}. Two appendices will discuss the extension to the multi-directional flow routing schemes, and outline considerations for implementations of the algorithms described herein.

In general, I will provide slightly more detail than one might when writing for readers with a background in scientific or parallel computing, in hopes of making the material accessible to those with a background in hydrology.

\paragraph*{A note on terminology.}
Whereas the term ``flow routing'' refers to the task of determining how much water flows where based on knowledge of how much rain fell elsewhere, the related problem of ``drainage area determination'' asks about the area that lies \textit{upstream} of a point. (Or perhaps more correctly: It asks for the area surrounding the nodes that drain through the area surrounding a given node.) The latter problem is a special case of the former: If one assumes a spatially constant rainfall of, say, one cubic meter of water per square meter of land, and finds by solving the flow routing problem that through a specific point $10^8$ cubic meters of water drain, then the drainage area of that point is clearly $10^8$ square meters (100 square kilometers). In other words, determining drainage area corresponds to the special case of flow routing with spatially constant rainfall.

\section{The parallel flow routing algorithm}
\label{sec:basic-algorithm}

In this section, let me first outline the sequential algorithm that is widely used (Section~\ref{sec:high-to-low}), along with the reformulation of flow routing as a system of linear equations (Section~\ref{sec:flow-routing-as-linear-systems}), and then the IDA algorithm first described by Richardson, Hill, and Perron in \cite{Richardson2014} (Section~\ref{sec:ida}). All three of these sections cover material that is well known, but that will combine to form the the new and much improved algorithm I will outline in Section~\ref{sec:better-algorithm} and then further optimize in Section~\ref{sec:optimized-algorithm}.

\subsection{Classical flow routing: The sequential ``high-to-low'' algorithm}
\label{sec:high-to-low}

\begin{figure}
    \centering
    \includegraphics[width=0.4\linewidth]{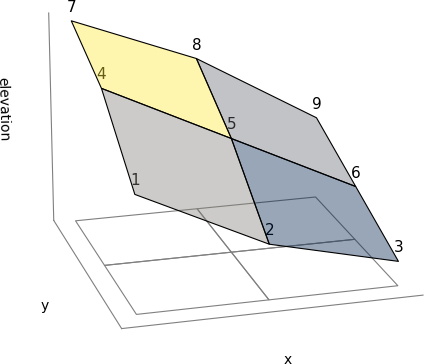}
    \hfill
    \includegraphics[width=0.4\linewidth]{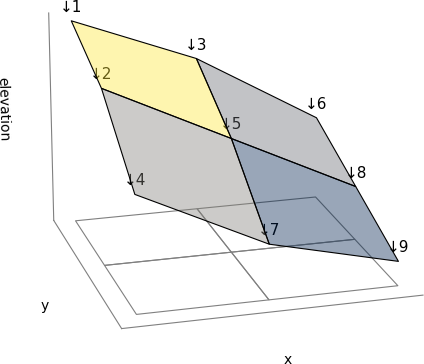}
    \caption{A simple $3\times 3$ digital elevation model. Left: Nodes are numbered left-to-right and bottom-to-top. Right: Nodes are numbered high-to-low, as indicated by the $\downarrow$ symbol next to each index.}
    \label{fig:3x3}
\end{figure}

In order to explain the classical ``high-to-low'' algorithm, consider the simple elevation model on a grid of $N=3\times 3=9$ points shown in Figure~\ref{fig:3x3}, and using the numbering of nodes of the left panel. The amount of water $w_i$ that flows through node $i$ is given by two contributions: First, the amount of rain $r_i$ that falls on the area around node $i$ (with this area defined, for example, by taking all points that are closer to node $i$ than any other node, which here corresponds to little rectangles around each node). And second, the amount of water any of the upstream nodes give to node $i$. For the second contribution, we need to determine for each node which other nodes it gives its water to; typically, in the common $D_4$ and $D_8$ schemes \cite{OCallaghanMark1984, Jenson1988}, one finds that one among the 4 immediate or 8 adjacent nodes that (i) is lower than node $i$ and (ii) to which the downhill gradient is the steepest.%
\footnote{Indeed, what matters for the discussion that follows is that \textit{each node gives water to exactly one neighbor} (or none, if at the boundary). The $D_4$ and $D_8$ schemes satisfy this requirement, but so would any other generalization on unstructured, triangular, or hexahedral meshes, or in schemes where water flows to a neighboring node that happens to not be the lowest one or not the steepest downhill. In the following, I will therefore reference the $D_4$ and $D_8$ schemes only in the sense of ``representatives'' of these kinds of schemes. There are of course other schemes in which the water from each node is partitioned among all of its neighbors in some way; I will comment on the generalization to such schemes in \ref{sec:Dinfty}.}

The step of computing which other node each node gives its water to is called ``local flow routing''. The output of local flow routing is a list of `(\textit{src}$\rightarrow$\textit{dst})' relationships that for the example in Fig.~\ref{fig:3x3} reads as follows:
\begin{align*}
&(1 \rightarrow 2), 
&&(2 \rightarrow 3), 
&&(3 \rightarrow \times), 
&&(4 \rightarrow 2), 
&&(5 \rightarrow 3), \\
&(6 \rightarrow 3), 
&&(7 \rightarrow 5), 
&&(8 \rightarrow 6), 
&&(9 \rightarrow 6).
&&
\end{align*}
Note that every node appears exactly once in the \textit{src} position. Node 3 is the lowest one and has no lower-lying neighbor to give its water to and so the `\textit{dst}' position is denoted by an ``invalid'' marker $\times$.%
\footnote{Nodes without lower-lying neighbors can be at the boundary where they would indicate the presence of a stream flowing out of the domain, or in the interior of the domain where they would represent local depressions into which water can flow but from which water can not escape without modifying the model. In flow routing, one avoids the presence of the second kind by pre-processing the digital elevation model in a ``depression-filling'' step. See also Section~\ref{sec:evaluation}.}

The key realization of the classical algorithm is that we can only determine how much water a node $i$ has if we know how much it gets from others, but that the only ones we need to care about are the \textit{higher-lying neighboring nodes}, which must consequently be processed earlier. In other words, if we process first the highest-lying node (which cannot have any inputs other than rain), then the second highest node (which either has no higher-lying inputs, or the one just processed), and then continue in this order, we will be able to process each node in a way that creates no depencies that are not already met by the time we touch a node. As a consequence, we have to order nodes from high to low and then process them in this order. Given the flow directions shown above, and sorting nodes from high to low, we can then write the amount of water available at every one of the 9 nodes of the example as follows:
\begin{align}
\label{eq:high-to-low-procedure}
\begin{split}
w_7 &= r_7, \\
w_4 &= r_4, \\
w_8 &= r_8, \\
w_1 &= r_1, \\
w_5 &= r_5 + w_7, \\
w_9 &= r_9, \\
w_2 &= r_2 + w_1 + w_4, \\
w_6 &= r_6 + w_8 + w_9, \\
w_3 &= r_3 + w_2 + w_5 + w_6.    
\end{split}
\end{align}
Here, the $r_i$ are all assumed to be known rainfall amounts, and the right hand side of the equation for $w_i$ contains as inputs all of the nodes $j$ that give water to node $i$ -- i.e., all of those `(\textit{src}$\rightarrow$\textit{dst})' relationships in which $i$ appears in the `\textit{dst}' position. Importantly, in each line, all of the $w_j$ that appear on the right hand side have already been computed by one of the \textit{previous} equations -- because we have sorted the list from high to low. In other words, we can solve for all 9 of the $w_i$ in turn, by running through the list of equations from top to bottom. We will thus call this algorithm the ``high-to-low'' algorithm. 

A simple piece of Python-like code that implements this algorithm -- building up the right hand sides of the equations as we go -- might look as follows:
\begin{codebox}[]{global flow routing}
\begin{lstlisting}

def global_flow_routing (local_flow_routing_list, r):
  w = r                                  # Start w vector with rainfall in r
  h2l_local_flow_routing = sorted(local_flow_routing_list, reverse=True)
  for src,dst in h2l_local_flow_routing: # Handle (src,dst) pairs high->low
    if dst != invalid:
      w(dst) = w(dst)+w(src)              # Add 'src's water to 'dst' node
  return w
\end{lstlisting}
\end{codebox}

It is difficult to imagine a more efficient (or more elegant!) algorithm for flow routing. There are two components to its run time: First, the cost of sorting all nodes from high to low, which for $N$ nodes can be achieved in a number of instructions that is proportional to $N\, \log_2(N)$ (we then say that sorting costs ${\cal O}(N \, \log_2(N))$ operations). Second, the actual traversal of the $N$ equations to compute each of the $w_i$ comes at a cost that is proportional to $N$. This is because there are $N$ equations and on the right hand side, there can be at most $4$ (for the $D_4$ scheme) or $8$ (for the $D_8$ scheme) terms, plus the one term for $r_i$, implying that the total cost is at most proportional to $5N$ or $9N$. Even for large models, $\log_2(N)$ is a relatively modest number -- for example, the largest model considered in Section~\ref{sec:evaluation} has $N=\num{1.88e9}$, and consequently $N\,\log_2(N)\approx 31 N$. In practice the two parts of the algorithmic cost are therefore of comparable size.

Various attempts have been made over the years to parallelize this algorithm; Section 2.4 of \cite{Richardson2014} provides a nice overview of approaches, and it is also worth mentioning \cite{Arge2003,Braun2013,Barnes2016} in this context. The key issue is that we traverse nodes from high to low, an inherently sequential order. One can try to split the work onto separate processes by ensuring that each process only works on nodes whose upstream nodes have either already been worked on, or are also dealt with by that process. Or one can keep track how many dependencies of a node are still unfilled, and keep a list of ``ready'' nodes for which all higher-lying neighbors have already been processed; the nodes in this list can then be worked on in parallel. These approaches have, however, largely not been very successful because (i) they do not scale to the large number of processes available on today's workstations and supercomputers, (ii) they introduce synchronization overhead that negates much of the benefit of parallel computations, and (iii) they typically required all processes to have access to the \textit{entire} data set in a shared-memory setting -- a major impediment at a time when data sets may have $N=10^{11}$ points and require terabytes of memory just to store.

As a consequence, one needs to approach flow routing in an entirely different way to make use of the large numbers of processor cores available today. The following section provides the conceptual basis for such an approach.

\subsection{Writing flow routing as a system of linear equations}
\label{sec:flow-routing-as-linear-systems}

The starting point for developing a parallel algorithm is the realization that flow routing can be written as a system of linear equations. As mentioned in the introduction, this connection has been made multiple times, apparently independently, in the literature and is not new. Considering the example of the previous section, let us bring all of the unknowns $w_i$ to the left side and sort by index:
\begin{align*}
w_1 				&= r_1 ,\\
w_2 - w_1 - w_4 		&= r_2, \\
w_3 - w_2 - w_5 - w_6 	        &= r_3, \\
w_4 				&= r_4, \\
w_5 - w_7 			&= r_5, \\
w_6 - w_8 - w_9 		&= r_6, \\
w_7 				&= r_7, \\
w_8 				&= r_8, \\
w_9 				&= r_9.
\end{align*}
The set of linear equations above can in turn be written in matrix-vector form as follows, where for clarity we have omitted the many zeros that fill the remainder of the matrix:
\begin{align}
\label{eq:linear-system}
\begin{pmatrix}
 1 &&&&&&&& \\
 -1 & 1 & & -1 &&&&& \\
 & -1 & 1 & & -1 & -1 &&& \\
 &&&1&&&&& \\
 &&&&1&&-1& \\
 &&&&&1&&-1&-1\\
 &&&&&&1&& \\
 &&&&&&&1& \\
 &&&&&&&&1
\end{pmatrix}
\begin{pmatrix}
    w_1 \\ 
    w_2 \\
    w_3 \\ 
    w_4 \\
    w_5 \\ 
    w_6 \\
    w_7 \\ 
    w_8 \\
    w_9
\end{pmatrix}
=
\begin{pmatrix}
    r_1 \\ 
    r_2 \\
    r_3 \\ 
    r_4 \\
    r_5 \\ 
    r_6 \\
    r_7 \\ 
    r_8 \\
    r_9
\end{pmatrix}.
\end{align}
In the following, we will abbreviate such linear systems as
\begin{align}
\label{eq:Aw=f}
    A \mathbf w = \mathbf r,
\end{align}
where $A$ is an $N\times N$ matrix, $\mathbf w$ is the vector of size $N$ that holds the water amounts of all nodes and that we would like to determine, and $\mathbf r$ is the known vector that for each of the $N$ node stores the local amount of rainfall.

At this point, we have formulated the flow routing problem in a way that allows us to utilize the universe of methods used to solve systems of linear equations. It is useful to recognize that the matrix is ``sparse'', i.e., that the vast majority of entries are zero. In fact, because each node in the $D_4$ and $D_8$ schemes gives its water to exactly one other node (or to none, if it is a local depression or a low point at the boundary of the domain), each column has at most two entries; each row $i$ contains a one on the diagonal and as many minus ones as the number of neighbors of $i$ that give their water to node $i$ -- i.e., at most 4 or 8 minus ones per row for the two schemes.

While it is \textit{elegant} to re-formulate the problem as a linear system, in itself this is not \textit{efficient}: In general, solving a linear system in $N$ variables takes a run time that grows proportional to $N^3$ (i.e., we say the algorithm is ${\cal O}(N^3)$). In the current case, one can make use of the sparsity of the matrix and other structural properties of the problem to substantially reduce the cost, but using a black-box linear system solver will almost certainly not come even close to the sequential algorithm outlined in Section~\ref{sec:high-to-low} that has a run time that is ${\cal O}(N \, \log_2(N)$. Instead, the value of the reformulation as a linear system will be that we can solve the problem \textit{in parallel} on multiple processes at once.

\subsection{The IDA algorithm: Solving systems of linear equations in parallel}
\label{sec:ida}

The innovation of \cite{Richardson2014} is to realize that there are algorithms that can solve linear systems such as \eqref{eq:linear-system}--\eqref{eq:Aw=f} in parallel, and that implementations of these algorithms are widely available in libraries such as those used in their paper. In general, these algorithms are iterative in nature, that is they consist of a loop in which the current best guess for the solution vector $\mathbf w$ is improved in each iteration. Iterative algorithms generally only require the computation of matrix-vector products with the matrix $A$ in \eqref{eq:Aw=f}, which can efficiently be done in parallel. In practice, all of these iterative algorithms require fewer operations when used in combination with a ``preconditioner'' -- a matrix $B$ that is used to transform the linear system $A\mathbf w=\mathbf r$ into the equivalent system $BA\mathbf w=B\mathbf r$ so that the matrix $BA$ has better properties than $A$ alone. In practice, it is not necessary to form the matrix product $BA$; instead, it is enough to form matrix-vector products $BA\mathbf x=B(A\mathbf x)$, i.e., multiply a vector $\mathbf x$ first by the matrix $A$ and then by $B$.

Based on this realization, \cite{Richardson2014} then proceeded to run tests on a wide range of linear solver algorithms and preconditioners to solve equations \eqref{eq:linear-system}--\eqref{eq:Aw=f} on multiple processes as well as on Graphics Processing Units (GPUs). How well each combination of solver and preconditioner worked is then summarized in a table in Figure 3 of \cite{Richardson2014}.
They find that the combination of a Richardson solver and the Euclid preconditioner performed the best. They then called this combination of writing flow routing as a linear system and solving it with Richardson+Euclid the ``Implicit Drainage Area'' (IDA) algorithm.%
\footnote{It is perhaps confusing that Richardson et al.~find that the \textit{Richardson solver} is the best method. I will be careful in the following when using the name to make clear whether I reference the authors of the paper, or the method. The method to solve linear systems was proposed by Lewis Fry Richardson in \cite{Richardson1911}, more than 100 years before the paper I reference frequently herein.}

\begin{figure}
    \centering
    \includegraphics[height=0.3\linewidth]{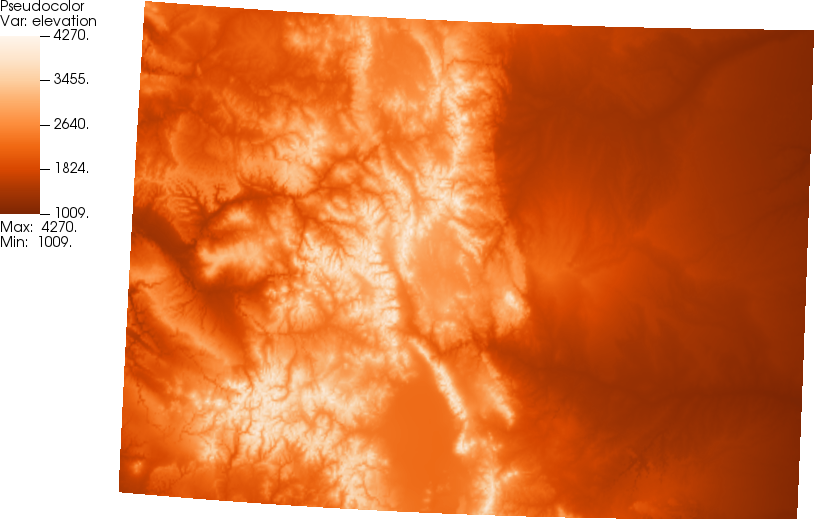}
    \hfill
    \includegraphics[height=0.3\linewidth]{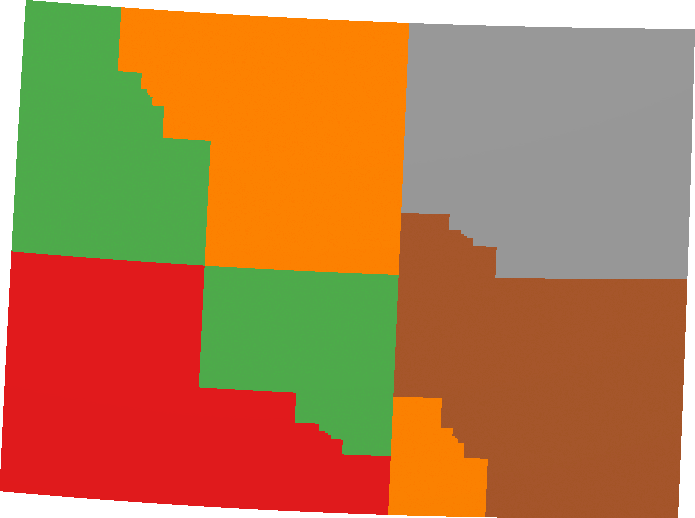}
    \caption{Left panel: A digital digital elevation model with $N=897\times 513=\num{460161}$ points of the state of Colorado, USA, the home state of the author. The model is shown projected onto the surface of Earth in three dimensions and clearly shows the Rocky Mountains in the western half of the state, and the valleys of the South Platte river existing at top right, of the Arkansas river existing at bottom right, and of the Colorado river exiting center left. Right panel: A partition of the node points among five processes, each indicated by a separate color. Note that partitions do not have to be neatly shaped or even to be contiguous.}
    \label{fig:colorado-dem}
\end{figure}

Richardson et al.~do not actually compare their algorithm when run on a single process against the very efficient high-to-low algorithm of Section~\ref{sec:high-to-low}. In fact, reading between lines and based on our own results in Section~\ref{sec:evaluation}, one might speculate that if a data set is small enough to be stored and run on a single process, the IDA algorithm is substantially slower than the high-to-low algorithm.
But this does not invalidate the approach: First, an algorithm that on a single process may be slower than the high-to-low one may still be much faster when run on sufficiently many processes, as we will see in Section~\ref{sec:evaluation}.
Moreover, iterative algorithms do not require that the matrix or solution vectors are stored on a single machine -- on a cluster of computers, each machine may only store a subset of the rows of the matrix, and the corresponding parts of the $\mathbf w$ and $\mathbf r$ vectors. As a consequence, one can deal with very large digital elevation models that cannot even be stored in memory on a single machine. In such cases, one would ``partition'' the domain of interest among the processes as shown in the right panel of Fig.~\ref{fig:colorado-dem}; each process in a parallel computation then ``owns'' only a subset of the model and will only need to query the DEM for the points it owns, plus one layer of ``ghost nodes'' around the ones it owns. There are widely used software libraries that can efficiently partition domains like the one shown into tens or hundreds of thousands of pieces.

\subsection{A better parallel algorithm for flow routing: Optimal preconditioners}
\label{sec:better-algorithm}

Some outcomes of the comparison of solvers and preconditioners in Figure 3 of \cite{Richardson2014} are not a surprise. That is because many of the solvers tested there -- for example, Conjugate Gradients (CG) -- require the matrix $A$ in \eqref{eq:Aw=f} to be symmetric and positive definite, but it is neither. Others of the methods tested offer flexibility in the choice of operations that are not necessary here (for example, F-GMRES over GMRES). Similarly, some of the preconditioners require the matrix to have properties it does not have in the current application. As a consequence, the finding that the Richardson solver works well is not surprising, though experts in the solution of linear systems would probably have chosen GMRES at first given its well understood performance on non-symmetric matrices.%
\footnote{Despite my initial preference for GMRES over the Richardson solver, my experiments mirror those of \cite{Richardson2014}: The two require the same number of iterations to reach the desired tolerance, but the Richardson solver takes only about 2/3 of the run time. This is because the internal operations of the GMRES solver -- but not those of the Richardson solver -- are quite expensive and, given the simplicity of the matrix $A$, contribute noticeably to the overall solver time. As a consequence, I agree with the choice of the Richardson solver.}
Euclid as a preconditioner \cite{Hysom99,Hysom01} is a more enigmatic -- and uncommon -- choice but turns out to be a good method because it is fundamentally an incomplete LU decomposition and, as will become clear throughout the rest of this section, computing matrix decompositions is quite easy if the matrix is triangular, a property the matrix $A$ herein has when looked at in the right way.

\begin{figure}
    \phantom{.}
    \hfill
    \centering
    \includegraphics[width=0.25\linewidth]{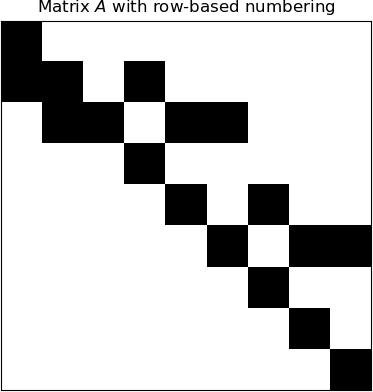}
    \hfill
    \includegraphics[width=0.25\linewidth]{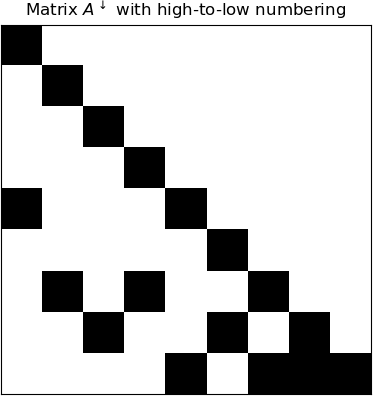}
    \hfill
    \phantom{.}
    \caption{A visualization of the nonzero entries of the matrix that appears in the reformulation of flow routing as a system of linear equations \eqref{eq:linear-system}--\eqref{eq:Aw=f}. Diagonal matrix entries have a value of $+1$, all other nonzero entries are $-1$. Left: The matrix $A$ we obtain if we enumerate nodes as in the left panel of Fig.~\ref{fig:3x3}, namely first in $x$- and then in $y$-direction. Right: The matrix $A^\downarrow$ we would have gotten instead if we had enumerated nodes high-to-low as shown in the right panel of Fig.~\ref{fig:3x3}. Note that $A^\downarrow$ is triangular.}
    \label{fig:matrix-3x3}
\end{figure}

Whereas the choice of the Richardson solver is reasonable, one can do substantially better with the preconditioner. Explaining how requires a short detour. In reformulating the $3\times 3$ flow routing problem shown in Fig.~\ref{fig:3x3} as a system of linear equations in matrix form \eqref{eq:linear-system}--\eqref{eq:Aw=f}, we derived the matrix $A$ whose nonzero entries can be visualized as in the left panel of Fig.~\ref{fig:matrix-3x3}. But we are not required to enumerate the nodes of the elevation model first in $x$- and then in $y$-direction. We could also have enumerated them high-to-low as shown in the right panel of Fig.~\ref{fig:3x3}, and in that case we would have obtained a matrix $A^\downarrow$ as shown in the right panel of Fig.~\ref{fig:matrix-3x3}. Depending on how we enumerate nodes, we can solve either $A\mathbf w=\mathbf r$ or $A^\downarrow\mathbf w^\downarrow=\mathbf r^\downarrow$ -- the two will lead to identical amounts of water at each node, it is just a matter what number we attach to each node. In fact, different enumerations of nodes are just permutations of indices, and if we denote the permutation by the matrix $P$, then $\mathbf w^\downarrow = P\mathbf w$, $\mathbf r^\downarrow = P\mathbf r$, and $ A^\downarrow = PAP^T$.

Despite the equivalence of the two linear systems, it is much easier to solve $A^\downarrow\mathbf w^\downarrow=\mathbf r^\downarrow$. This is because the matrix $A^\downarrow$ is \textit{(lower) triangular}, i.e., all entries above the diagonal are zero and all nonzeros are on or below the diagonal. Linear systems with triangular matrices are particularly easy to solve because we can use the first row of $A^\downarrow$ and the first entry of the right hand side, $r^\downarrow_1$, to solve for $w^\downarrow_1$. With $w^\downarrow_1$ known, we can use the second row of $A^\downarrow$ to solve for $w^\downarrow_2$, knowing that the entries $A^\downarrow_{2,3}\cdots A^\downarrow_{2,N}$ are all zero and so none of the water variables not yet computed appear in this equation. We can continue row by row to compute $\mathbf w^\downarrow$; this step-by-step solve is called ``forward substitution''. Of course, it is no surprise that solving the renumbered linear system is easy: We have numbered nodes high-to-low, and we know that we can easily solve for all water degrees of freedom in this order: This is what the high-to-low algorithm of Section~\ref{sec:high-to-low} does, after all! In other words, the simplicity of the high-to-low algorithm corresponds to the fact that, if we enumerate nodes smartly, the matrix we need to solve with in the formulation as a linear system is triangular and easy to solve with.

\begin{figure}
    \phantom{.}
    \hfill
    \centering
    \includegraphics[height=0.3\linewidth]{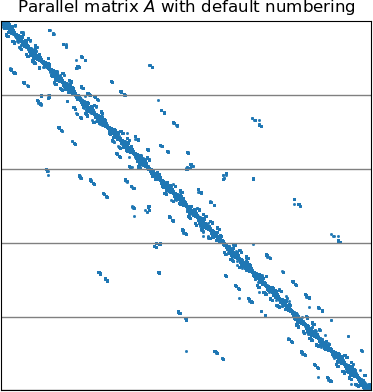}
    \hfill
    \includegraphics[height=0.3\linewidth]{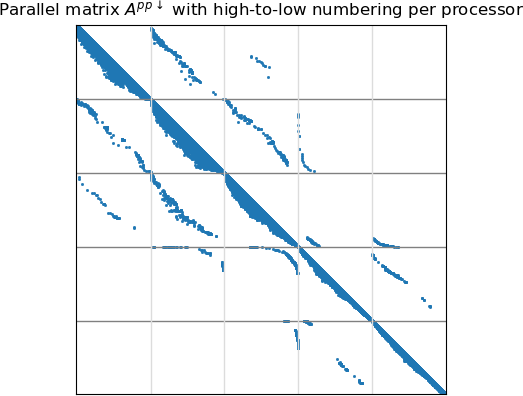}
    \hfill
    \phantom{.}
    \caption{The matrices $A$ (left) and $A^{pp\downarrow}$ (right) that correspond to the Colorado elevation model of  Fig.~\ref{fig:colorado-dem} with $N=\num{460161}$ points, partitioned onto $P=5$ processes. Horizontal gray horizontal lines indicate the blocks of rows stored by each of the processes. In the right panel, vertical light gray lines also indicate the corresponding partitioning of columns, though this has no effect on how matrices are built or stored.}
    \label{fig:matrix-dem}
\end{figure}

Let us now consider how this observation can be used in parallel computations. The basic assumptions we will have to live with are that (i) we have partitioned the domain among the $P$ processes, such as in the situation shown in the right panel of Fig.~\ref{fig:colorado-dem} with $P=5$, where each process ``owns'' the nodes that are part of its partition; and (ii) that each process $p=1,\ldots,P$ is responsible for a contiguous block of $N_p$ node indices that it can otherwise freely assign to the nodes it owns. In other words, process 1 will own indices $1\ldots N_1$ that it can assign to its $N_1$ nodes; process 2 will own indices $N_1+1\ldots N_1+N_2$ that it can assign to its $N_2$ nodes; and so on for a total of $N_1+N_2+\cdots+N_P=N$ indices. The second assumption ensures the practical consideration that each process should only have to store contiguous sets of rows of the matrix $A$ and contiguous sets of entries of vectors $\mathbf w, \mathbf r$. The left panel of Fig.~\ref{fig:matrix-dem} shows the matrix $A$ that corresponds to the Colorado elevation model in Fig.~\ref{fig:colorado-dem} when partitioned onto $P=5$ processes, using the default numbering of nodes provided by the \dealii{} library used for the computations here \cite{arndt2019dealii,Arndt2025}.

Because of the requirement that each parallel process owns contiguous ranges of nodes, we cannot \textit{globally} renumber high-to-low. But each process can independently and in parallel enumerate its own nodes high-to-low. We will call the resulting matrix $A^{pp\downarrow}$ (``per-process high-to-low''); for the Colorado model, the nonzeros of this matrix are shown in the right panel of Fig.~\ref{fig:matrix-dem}. This matrix is not triangular because a node might receive water from a higher-lying node that is owned by one of the later processes; that higher-lying node therefore has a larger index (unlike in the one-process case where it would have received a lower index) and so leads to a matrix entry to the right of the diagonal.

Coming back to the question of preconditioning the iterative solution of the linear system $A\mathbf w=\mathbf r$ (or $A^{pp\downarrow}\mathbf w^{pp\downarrow}=\mathbf r^{pp\downarrow}$), let us first recall that a preconditioner is a matrix $B$ that transforms the linear system to $BA\mathbf w=B\mathbf r$. The best preconditioner would be $B= A^{-1}$ given that that would convert the linear system into $\mathbf w=B\mathbf r$. In practice, of course, computing the matrix $A^{-1}$ itself cannot be done efficiently, and so one has to balance between $B$ being close to $A^{-1}$ and $B$ also being efficiently computable. \cite{Richardson2014} use the Euclid algorithm that chooses $B$ to be an incomplete LU decomposition of $A$ to approximate $A^{-1}$, whereas herein we will make use of the fact that we have a very efficient \textit{exact} solver, namely the high-to-low procedure. Second, all iterative solvers do not actually need the matrix $B$ entry-by-entry; rather, they just need the ability to compute matrix-vector products $\mathbf y=B\mathbf x$ given a vector $\mathbf x$.

In the current context, we note that while $A^{pp\downarrow}$ is not triangular, the \textit{diagonal blocks} of $A^{pp\downarrow}$ are triangular (see also the right panel of Fig.~\ref{fig:matrix-dem}) because the nodes of each process are individually numbered high-to-low. Let us then write $A^{pp\downarrow}$ as the sum of a matrix $D^{pp\downarrow}$ composed of the diagonal blocks, and a matrix $R^{pp\downarrow}$ that contains the rest, that is
\begin{align}
    \label{eq:matrix-sum}
    \underbrace{
      \includegraphics[height=2cm,valign=m]{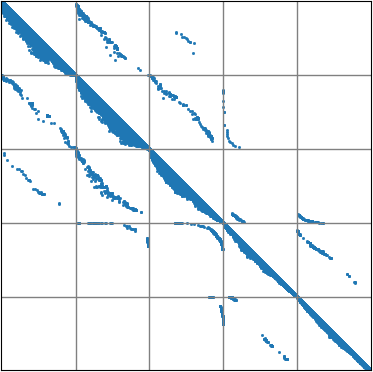}
    }_{A^{pp\downarrow}}
    =
    \underbrace{
      \includegraphics[height=2cm,valign=m]{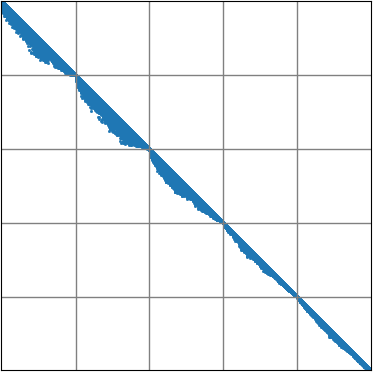}
    }_{D^{pp\downarrow}}
    +
    \underbrace{
      \includegraphics[height=2cm,valign=m]{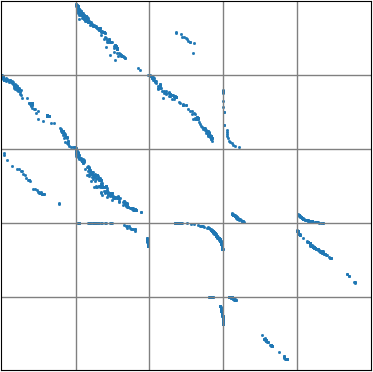}
    }_{R^{pp\downarrow}}.
\end{align}
We will then choose our preconditioner $B=(D^{pp\downarrow})^{-1}$. Computing matrix-vector products $\mathbf y=B\mathbf x=(D^{pp\downarrow})^{-1}\mathbf x$ is then equivalent to solving the linear system $D^{pp\downarrow}\mathbf y=\mathbf x$, which is easy because $D^{pp\downarrow}$ is triangular. Moreover, because it only consists of the diagonal blocks, \textit{each parallel process can solve its linear system independently of that of the other processes}, and in parallel. As discussed above, solving these triangular systems corresponds to applying the high-to-low procedure of Section~\ref{sec:high-to-low} on the nodes that are owned by a process.

In summary, our solution approach will then be the following:
\begin{enumerate}
    \item We use the Richardson method as the iterative solver.
    \item We use $B=(D^{pp\downarrow})^{-1}$ as the preconditioner.
\end{enumerate}
It is worth outlining our expectations how this method works. First, if we compute with only a single process (i.e., we are performing a ``sequential'' computation), then all matrices consist of only a single block and so $A^{pp\downarrow}=A^\downarrow$, $D^{pp\downarrow}=A^{pp\downarrow}$, $R^{pp\downarrow}=0$, and finally $B=(A^\downarrow)^{-1}$. This means that the preconditioned linear system is $BA^\downarrow \mathbf w = I \mathbf w = B\mathbf r$. The Richardson solver converges in one step whenever the matrix in the linear system is the identity matrix $I$, and so the cost of solving the linear system is one multiplication by $B$ -- i.e., equivalent to one high-to-low solve. This is pleasing: In the sequential case, our method is as efficient as the best algorithm we have for flow routing.

In contrast, for $P>1$ parallel processes, $B\neq (A^\downarrow)^{-1}$ and we should expect the Richardson solver to take more than one iteration. In particular, as we add more and more processes, we should expect that $R^{pp\downarrow}$ consists of a larger and larger portion of the matrix $A^{pp\downarrow}$, and consequently that $B=(D^{pp\downarrow})^{-1}$ becomes a worse and worse approximation of $(A^{pp\downarrow})^{-1}$. As a consequence, for a problem of fixed size, we should expect that the number of Richardson iterations necessary to solve for $\mathbf w^{pp\downarrow}$ becomes larger and larger. We will show that this expectation is indeed correct in Section~\ref{sec:evaluation-iterations} and provide some theoretical insight why this is so in Section~\ref{sec:richardson}.

Let me end this section by stating that the method described above can be interpreted in the spirit of domain decomposition (DD) approaches \cite{QV99} in which the subdomain solvers (here: the high-to-low solver for each partition) are exact. This has no bearing for the remainder of the paper, but could in principle be used to develop better methods if one could come up with interface transmission conditions that allow for building solvers that work on the skeleton of the domain decomposition.

\section{Further optimization of the algorithm}
\label{sec:optimized-algorithm}

Having outlined the basic ideas of a better parallel flow routing algorithm above, let us discuss in this section a number of improvements one can make to algorithm. Some of these are algorithmic in nature (specifically, the ones outlined in Sections~\ref{sec:optimization-matrix-free} and \ref{sec:optimization-unsorted}) whereas others modify the mathematical formulation (Sections~\ref{sec:optimization-combined-I+X} and \ref{sec:optimization-I+X-times-I-X}).

\subsection{Optimization 1: A ``matrix-free'' algorithm}
\label{sec:optimization-matrix-free}

As discussed in Section~\ref{sec:better-algorithm}, the key insight of this paper is that using the per-process high-to-low operation is the optimal choice as a preconditioner for the formulation of flow routing as a linear problem that can be solved with parallel iterative solvers. Like all other iterative methods, the Richardson solver does not actually require access to specific matrix entries; instead, all that it needs is the ability to multiply the matrix by a vector, and similarly to apply the preconditioner to a vector. In other words, these methods do not actually care about how the matrix and preconditioner are stored, but have a \textit{functional} perspective of these operators in that they only need to know how the operators can be \textit{applied} to a vector.

A common approach optimizing the time to solution in such cases is to not actually store the matrix in one of the common sparse matrix storage formats (say, in Compressed Row Storage/CSR format), but to work in a ``matrix-free'' way in which one implements the \textit{action} of the matrix on a vector. In the current context, this is easily possible.

To understand how one would do that, note that in building the linear system, we needed to first find for each node which other node (if any) it gives water to (i.e., perform ``local flow routing''). As discussed in Section~\ref{sec:high-to-low}, the result of local flow routing is a list of pairs (\textit{src}$\rightarrow$\textit{dst}) that for each node `\textit{src}' says where water is routed to.

The key insight is that each (\textit{src}$\rightarrow$\textit{dst}) pair corresponds to one \textit{column} of the matrix: As can be seen in the derivations in Section~\ref{sec:flow-routing-as-linear-systems}, each node `\textit{src}' appears with a $+1$ on the diagonal of the matrix in column `\textit{src}', and with a $-1$ in row `\textit{dst}' of column `\textit{src}', the latter indicating that `\textit{src}' is one (of possibly several) sources in the equation for `\textit{dst}' represented by row `\textit{dst}' in the matrix. Because in the $D_4$ and $D_8$ schemes, each node only gives water to one other node, these are the only two entries in column `\textit{src}' of the matrix. (In the $D_\infty$ and related schemes, each node can give water to several others; we discuss the extension of the approach outlined here in \ref{sec:Dinfty}.) 

Using this column-oriented approach, we can implement the operation $y=Ax$ that forms the core of what iterative solvers do with the following algorithm (again, written in Python-like syntax):
\begin{codebox}[]{matrix-vector product with $A$}    
\begin{lstlisting}
def apply_A (x):
  y = 0                                     # Set output vector to zero
  for src,dst in local_flow_routing_list:   # Loop over (src->dst) list
    y(src) = y(src) + A(src,src) * x(src)   # Note: A(src,src) equals +1
    if dst != invalid:
      y(dst) = y(dst) + A(dst,src) * x(src) # Note: A(dst,src) equals -1
  return y
\end{lstlisting}
\end{codebox}
Here, we assume that the \texttt{local\_flow\_routing\_list} variable contains the (\textit{src}$\rightarrow$\textit{dst}) pairs.
As mentioned in the comments of the listing, one can further optimize the algorithm by noting that the matrix entries are known and need not actually be stored. The following algorithm is therefore equivalent:
\begin{codebox}[]{matrix-vector product with $A$; optimized}
\begin{lstlisting}
def apply_A (x):
  y = x                                    # Set y = diagonal of A times x
  for src,dst in local_flow_routing_list:  # Loop over (src->dst) list
    if dst != invalid:
      y(dst) = y(dst) - x(src)             # Handle the -1 in column 'src'
  return y
\end{lstlisting}
\end{codebox}
In the parallel context, each process is only responsible for a subset of rows of the matrix (and entries in the output vector); the algorithm is easily adjusted to that case by limiting \texttt{local\_flow\_routing\_list} to only contain information about nodes owned by the current process. Indeed, these are the only nodes about which each process has sufficient knowledge to compute local flow routing information anyway.

Similarly, one can write the preconditioner as a function without actually storing the matrix in question. Recall that the preconditioner $B=(D^{pp\downarrow})^{-1}$ is the high-to-low operation on each parallel process's diagonal block; these diagonal blocks of $A$ are triangular. In other words, if we denote a process's diagonal block by $L$ (to indicate that it is lower-triangular), then we want to compute $\mathbf y=L^{-1}\mathbf x$, which is equivalent to solving the linear system $L\mathbf y=\mathbf x$ for $\mathbf y$ where $\mathbf x$ is a given vector. Let us illustrate this with a small $3\times 3$ case where we need to solve
\begin{align*}
  \begin{pmatrix} 
  L_{11} & 0 & 0 \\
  L_{21} & L_{22} & 0 \\
  L_{31} & L_{32} & L_{33}
  \end{pmatrix} 
  \begin{pmatrix} 
  y_1 \\ y_2 \\ y_3
  \end{pmatrix} 
  =
  \begin{pmatrix} 
  x_1 \\ x_2 \\ x_3
  \end{pmatrix}.
\end{align*}
Here, we first solve $y_1 = \frac{x_1}{L_{11}}$, after which we can use the now known $y_1$ value in the second equation to compute $y_2 = \frac{x_2-L_{21}y_1}{L_{22}}$, and finally $y_3 = \frac{x_3-L_{31}y_1-L_{32}y_2}{L_{33}}$. On matrices of arbitrary size $n$, one would iterate over $k=1\ldots n$ and compute 
\begin{align*}
    y_k = \frac{x_k-\sum_{l=1}^{k-1}L_{kl}y_l}{L_{kk}}.
\end{align*}
While correct, this order of operations requires us to work row-based when computing the sum $-\sum_{l=1}^{k-1}L_{kl}y_l$ because it accesses the elements of row $k$ of the matrix. However, local flow routing has given us column-based information. It is not difficult to modify the algorithm for this: We just need to recognize that based on the local flow routing information, only one entry below $L_{kk}$ can be nonzero, say $L_{pk}$ in row $p>k$. So, if we introduce a temporary vector $\mathbf z$, we can add the term $-L_{pk}y_k$ that will appear in the equation for $y_p$ to $z_p$. 

Applied to the situation where we want to solve with the triangular diagonal blocks of the matrix $A$ instead of a generic $L$ then results in the following code:
\begin{codebox}[]{triangular solve}
\begin{lstlisting}
def solve_with_diagonal_blocks (x):
  y = 0                                      # Set output to zero
  z = 0                                      # Set temp vector to zero
  for src,dst in local_flow_routing_list:    # Loop over (src->dst) list
    y(src) = (x(src)+z(src)) / A(src,src)    # Solve for the 'src'th row
    if dst != invalid:
      z(dst) = z(dst) - A(dst,src) * y(src)  # Update z(dst)
  return y
\end{lstlisting}
\end{codebox}
This function assumes that the \texttt{local\_flow\_routing\_list} is sorted in ascending order of \textit{src}, which because we have ordered nodes high-to-low means that we process nodes in high-to-low order.
 
As before, we have excluded the case where a node has no downstream neighbor. Again, we can optimize the code by (i) noting that we know that the values stored in the matrix $A$ are either $+1$ or $-1$; and (ii) noting that based on the high-to-low sorting of unknowns, we necessarily have \textit{dst}$>$\textit{src} which implies that we write into elements of $\mathbf z$ further down than the one we are currently working on in $\mathbf y$ -- i.e., we can use the elements of $\mathbf y$ as the necessary temporary storage! In consequence, the following optimized algorithm is equivalent to the one above:
\begin{codebox}[]{triangular solve; optimized}
\begin{lstlisting}
def solve_with_diagonal_blocks (x):
  y = 0                                      # Set output to zero
  for src,dst in local_flow_routing_list:    # Loop over (src->dst) list
    y(src) = x(src)+y(src)                   # Solve for the 'src'th row
    if dst != invalid:
      y(dst) = y(dst) + y(src)               # Update temporary storage in
                                             # as-yet unused elements of y
                                             # below y(src).
  return y
\end{lstlisting}
\end{codebox}
As before, it is not difficult to restrict this operation to only that part of the matrix that is actually owned by a process -- i.e., the diagonal blocks in the right panel of Figure~\ref{fig:matrix-dem}. The preconditioner discussed in Section~\ref{sec:better-algorithm} then applies these diagonal blocks to the sets of nodes stored by each process in a parallel computation.

The end result of these considerations is that we have two functions that implement the \textit{action} of the matrix and preconditioner -- which is all that iterative solvers require. It is not necessary to actually allocate memory to store the sparsity pattern and the entries of the matrix $A$. As we will see in Section~\ref{sec:evaluation}, this ``matrix-free'' scheme may be faster or slower than the matrix-based approach, depending on the characteristics of the machine one runs on. However, the matrix-free framework will allow us to consider additional optimizations in Sections~\ref{sec:optimization-combined-I+X} and \ref{sec:optimization-I+X-times-I-X} that are better or comparable to the matrix-based approach on all platforms tested.

\subsection{Optimization 2: Not enumerating degrees of freedom high-to-low}
\label{sec:optimization-unsorted}

The description of the previous sections assumed that nodes are numbered high to low, and that the \texttt{local\_flow\_routing\_list} variable that appears in the listings above contains the (\textit{src}$\rightarrow$\textit{dst}) pairs sorted in ascending order of \textit{src}. In that case, working through the pairs of the list guarantees that we are traversing the domain high-to-low.

On the other hand, in typical DEM models, nodes are ordered in some other order -- for example, using a lexicographic latitude-longitude scheme; or in computational models that work on something that is not a rectangular grid, the order may not have any geometric basis. \textit{Numbering} nodes in a high-to-low order therefore requires attaching another index to each node, and storing the high-to-low order in it. However, all of that is not necessary: In the code snippets shown above, what is required is that we \textit{process} the \texttt{local\_flow\_routing\_list} in high-to-low order; nodes need not be sorted in that order, and this assumption was only used in previous sections to illustrate that \textit{if one did it that way}, the matrix would have triangular structure.

As a consequence, we will no longer assume that node indices are sorted by elevation in the sections below, and correspondingly drop the $pp\downarrow$ superscript on all symbols from here on out. The only assumption we make is that the \texttt{local\_flow\_routing\_list} variable is sorted in high-to-low order of the \textit{src} part of the pairs, guaranteeing that we only ever process nodes whenever all higher neighbors have been processed.

\subsection{Optimization 3: Combining matrix and preconditioner -- the ``$I+X$'' scheme}
\label{sec:optimization-combined-I+X}

Conceptually, using a (left) preconditioner $B$ when solving a linear system $A\mathbf w=\mathbf r$ can be understood as solving the linear system
\begin{align}
\label{eq:BAw=Br}
  BA\mathbf w = B\mathbf r
\end{align}
instead. This is mathematically not entirely equivalent to what some iterative algorithms -- say, GMRES -- \textit{actually} do with preconditioners, but close enough for our purposes. In the current context, recall that $A$ has a structure as shown in Fig.~\ref{fig:matrix-dem} and that in \eqref{eq:matrix-sum} we have written it as
$A=D+R$ where $D$ corresponds to the diagonal blocks and $R$ to the off-diagonal blocks. $R$ has very few entries: in a parallel setting, each entry in $R$ corresponds to an off-process node that gives water to an on-process node. If we are working on only one process, then of course we have $A=D$ and $R=0$.

Next, consider that the preconditioner introduced in Section~\ref{sec:better-algorithm} is $B=D^{-1}$ and so $BA=D^{-1}(D+R)=I+D^{-1}R$ where $I$ is the identity matrix. That is, the preconditioned linear system $BAx=Bb$ is really
\begin{align}
\label{eq:I+X}
    (I+X) \mathbf w = D^{-1}\mathbf r \qquad\qquad \text{with $X=D^{-1}R$}.
\end{align}
Instead of the preconditioned system \eqref{eq:BAw=Br}, we might as well solve \eqref{eq:I+X}, obtaining the same solution as before. At first, one might think that for this it is necessary to form the matrix $X$, but as outlined in Section~\ref{sec:optimization-matrix-free}, this is not actually necessary: All we have to do is provide a procedure that implements the \textit{action} of $I+X$. This is not difficult, using the following function:
\begin{codebox}[]{multiply by $I+X$}
\begin{lstlisting}
def apply_IplusX (x):
  tmp = apply_R(x)                         # compute tmp=R*x
  y = solve_with_diagonal_blocks(tmp)      # compute y=D^{-1}*tmp
  y = y + x                                # add I*x=x
  return y
\end{lstlisting}
\end{codebox}
Here, \texttt{apply\_R} is a function similar to \texttt{apply\_A} discussed in Section~\ref{sec:optimization-matrix-free}, just restricted to matrix entries in off-diagonal blocks of $A$.
The reformulation in this way -- while otherwise entirely equivalent -- has two advantages: First, compared to the original formulation where we needed to multiply by $A$ and solve with triangulation systems $D$ (for the preconditioner) in each solver iteration, we now only need to multiply by $R$ and solve with $D$. But $R$ has many fewer entries than $A$, so multiplying with $R$ is substantially cheaper than with $A$. In the extreme case, where we work with only one process, we have $A=D$ and $R=0$, so $X=0$. Second, we can understand why the scheme described in Section~\ref{sec:better-algorithm} converges in relatively few iterations: We should think of $X$ as a ``small'' correction, and so we are solving a linear system $(I+X)\mathbf w=D^{-1}\mathbf r$ with a matrix that is not substantially different from the identity matrix for which \textit{any} iterative scheme should hopefully converge quickly!

\subsection{Optimization 4: A preconditioner for the combined system -- the ``$(I-X)(I+X)$'' scheme}
\label{sec:optimization-I+X-times-I-X}

If we solve the reformulated linear system $(I+X)\mathbf w = D^{-1}\mathbf r$ of the previous section with an iterative solver, then it would be nice if we could precondition this linear system in some way. The question is what a good preconditioner for the matrix $I+X$ might be, i.e., what matrix $Y\approx (I+X)^{-1}$.

The perspective of the previous section, namely that $X$ is a ``small'' correction to the identity matrix $I$ suggests an answer to this question: Based on the identity $\frac{1}{1+x}=(1+x)^{-1}=1-x+x^2-x^3+x^4-\cdots$, true for all $x$ with $|x|<1$ and with an obvious extension to matrices, one can approximate $(I+X)^{-1} \approx I-X =: Y$. Indeed, the preconditioned matrix $Y(I+X)=(I-X)(I+X)=I-X^2$ should be seen as being even closer to the identity matrix than the one in the previous section, and solving the preconditioned linear system
\begin{align}
    (I-X)(I+X)\mathbf w=(I-X)D^{-1}\mathbf r
\end{align}
should be possible in even fewer iterations than before.

Indeed, as we will show below, this is so: The number of iterations using the $I-X$ preconditioner compared to using no preconditioner at all on the $I+X$ matrix is reduced by a factor of about 2. On the other hand, applying $I-X$ costs as much as applying $I+X$, and so each iteration is twice as expensive, and in balance, the run times of the schemes using/not using $I-X$ as a preconditioner are comparable.

One could presumably further reduce the number of iterations by using $Y=I-X+X^2=I-X(I-X)$ or even higher order schemes that require several applications of $D^{-1}$ per preconditioning step. We leave this for future experimentation.

\section{Numerical evaluation of the parallel algorithm}
\label{sec:evaluation}

In the following, then, let us evaluate the four schemes for solving flow routing in parallel discussed in previous sections:
\begin{enumerate}
    \item The matrix-based scheme in which we build $A^{pp\downarrow}$ and precondition with $(D^{pp\downarrow})^{-1}$; see Section~\ref{sec:better-algorithm}.
    \item The matrix-free scheme with the same $A^{pp\downarrow}$ and $(D^{pp\downarrow})^{-1}$, but in which instead of storing the matrix, we only implement their action; see Section~\ref{sec:optimization-matrix-free}.
    \item The ``I+X'' scheme in which we combine matrix and preconditioner; see Section~\ref{sec:optimization-combined-I+X}.
    \item The ``(I-X)(I+X)'' scheme in which we  scheme in which we precondition the combined $I+X$ matrix by $I-X$; see Section~\ref{sec:optimization-I+X-times-I-X}.
\end{enumerate}
All four methods have been implemented in a code that is based on the widely used open source library \dealii \cite{arndt2019dealii,Arndt2025} that provides the infrastructure for handling and partitioning meshes in parallel, as well as many of the linear algebra building blocks used here (such as the Richardson solver). We build the parallel algorithms on the framework described in \cite{BBHK11}.
To facilitate easier dissemination of the methods described herein, the code itself along with extensive documentation is freely available as part of \dealii's ``code gallery'', a collection of contributed codes that typically accompany published papers; this code is currently available as a pull request at \url{https://github.com/dealii/code-gallery/pull/238} and I will make sure that it gets merged before the paper is published.

The program computes the water flow rate for the elevation model of Colorado -- the home state of the author -- shown in the left panel of Fig.~\ref{fig:colorado-dem}, assuming a spatially constant amount of rain of $r=\SI{375}{liter \per m^2 \per year}$, the approximately correct rainfall rate of Colorado's eastern plains. The digital elevation model (DEM) was obtained from OpenTopography and is based on the Copernicus Global Digital Elevation Models 90 meters, downsampled to 1800 meters, and then interpolated onto whatever resolution necessary for our numerical experiments. The data set was collected by the European Space Agency (ESA) through the Copernicus Programme, and in particular using data from the German TanDEM-X satellite. Specifics of this data set can be found at
\url{https://portal.opentopography.org/datasetMetadata?otCollectionID=OT.032021.4326.1}, see also \cite{CopernicusData}.

Based on this elevation model, we then create test cases at three different resolutions, resulting in meshes with $N=\num{3585}\times \num{2049}\approx \num{7.34e6}$, $N=\num{14337}\times \num{8193}\approx \num{117e6}$, and $N=\num{57345}\times \num{32769}\approx \num{1.88e9}$ node points,%
\footnote{These are all uniform subdivions by factors of 2 of an original $7\times 4$ cell mesh, given Colorado's extent of 7 by 4 degrees.} corresponding to surface resolutions of approximately \SI{160}{m}, \SI{40}{m}, and \SI{10}{m}, respectively. These models, which were originally given as elevation over a longitude-latitude grid, are then projected onto the three-dimensional surface of the earth to ensure that areas are correctly computed.
Before using the three meshes as the basis of evaluating the algorithms described herein, the program fills in depressions in the DEM via the priority-flood algorithm \cite{Barnes2014} to ensure that all water reaches the boundary of the domain.

\begin{figure}
    \centering
    \includegraphics[width=0.8\linewidth]{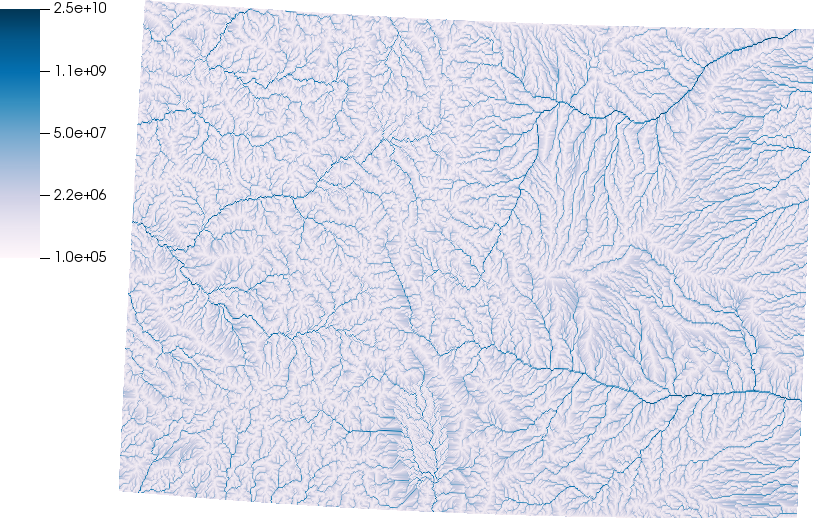}
    \caption{Computed water flow in \si{m^3 \per year} based on the digital elevation model shown in Fig.~\ref{fig:colorado-dem}. Colorado's three large rivers -- the South Platte, Arkansas, and Colorado rivers -- are clearly visible.}
    \label{fig:water}
\end{figure}

As with every software, it is useful to verify that what it computes is correct. Fig.~\ref{fig:water} shows the computed water flow rate, illustrating that at least qualitatively the computed solution is correct. To verify this also quantitatively, we note that the area covered by the model is $a=\SI{2.69016e+11}{m^2}=\SI{2.69016e+5}{km^2}$.%
\footnote{The \textit{true} area of Colorado slightly deviates from this: For the computations herein, we have chosen a domain that is delineated by 109 and 102 degrees West, and 37 and 41 North, as originally intended by the Colorado Organic Act of 1861 that created the territory. The \textit{actual} boundaries of the state, however, are defined by survey markers laid down in the 1860s and 1870s and deviate from the ideal boundary, sometimes by hundreds of meters, due to surveying errors.}
At the given rainfall rate, and since our model does not have a mechanism to remove water otherwise, we therefore expect the streams that cross the boundary to carry $r\cdot a=\SI{1.0088e11}{m^3 \per year}$ of water. Indeed, the computations show that the computed outflow rate equals this number to better than a relative error of $10^{-11}$. Moreover, the program checks that the solutions computed with the four different methods agree; they are in fact all equal to a relative difference of less than $10^{-6}$ in all experiments. Taken together, these are good indications for the correctness of the program.

In the following, let us then evaluate the four algorithms by two metrics: (i) The number of Richardson iterations it takes to compute the solution of the routing problem, in Section~\ref{sec:evaluation-iterations}; and (ii) the run time to solve the problem, in Section~\ref{sec:evaluation-runtime}. In both cases, the interesting question is how these numbers depend on the problem size and the number of processes used in the parallel computations.

\subsection{Number of iterations for the algorithms}
\label{sec:evaluation-iterations}

Let us first evaluate how many iterations are needed to solve the flow routing problem with the Richardson solver. This number is relevant because each iteration for the first three methods costs approximately as much in \textit{total} CPU cost as one high-to-low sweep in the sequential algorithm.%
\footnote{This is not quite correct. A high-to-low sweep routes water from all nodes to their respective neighbors, whereas the majority of the cost of one Richardson iteration is the application of the $D^{-1}$ preconditioner that only routes water between nodes \textit{within} each partition. In other words, it does marginally less work because it neglects routing across partition interfaces.}
However, this cost is now spread across multiple processes running in parallel. In other words, if the algorithm requires $n$ iterations to converge on $P$ processes, and assuming that communication between processes is negligible, we would expect a computation to take $\frac{n}{P}$ as much time as one high-to-low sweep on a single process. In order to see a speed-up in terms of the elapsed time required to solve the problem (the ``wall time'', referring to the clock hanging on the wall), we require that $n<P$.

\begin{figure}
    \centering
    \includegraphics[width=0.98\linewidth]{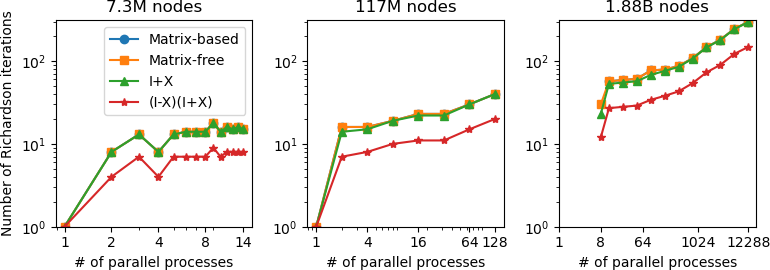}
    \caption{Number of Richardson iterations required to solve the flow routing problem. The three panels show results for digital elevation models of the Colorado topography with \num{7.3e6}, \num{117e6}, and \num{1.88e9} node points. The first three formulations result in nearly identical numbers of iterations. For the rightmost panel, computations are too large to be run with $P<12$ (matrix-based) or $P<8$ processes (for the three matrix-free variants).}
    \label{fig:iterations}
\end{figure}

Fig.~\ref{fig:iterations} then shows the number of iterations required for solving the flow routing problem using the Colorado DEM and the three meshes of different resolutions. The plots illustrate that all methods -- as expected -- solve the problem in one iteration when using a single parallel process; this must be so because in that case, the per-process high-to-low preconditioner $D^{-1}$ is the exact inverse of the matrix $A$, leaving us to solve the trivial linear system $I\mathbf w = D^{-1}\mathbf r$ for $\mathbf w$. On the other hand, the more parallel processes we use, the larger the number of iterations becomes. This should not come as a surprise: The preconditioner performs a high-to-low sweep over the nodes of each partition of the mesh only, and so gets the water flow rate at every node of the domain correct assuming that the inputs are correct. But the inputs include nodes on neighboring partitions that are not yet correct. Intuitively, one would expect that in the worst case, one needs to do as many iterations as the largest number of partitions a single stream crosses (see also the discussion in \ref{sec:implementation}). Algorithms like GMRES -- but not Richardson -- perform \textit{global} corrections that might alleviate this issue, but as mentioned above are in practice slower than using the Richardson solver.

A separate observation is that Methods 1 through 3 all use almost exactly the same number of iterations. This is not a surprise either, given that they are all reformulations of each other, but using the same underlying structure. Method 4 (the $I-X$ preconditioner applied to the $I+X$ matrix) requires pretty much exactly half as many iterations; that, too, is not overly surprising given that it is about twice as expensive, requiring two instead of one application of $D^{-1}$ per iteration.

Regardless of the growth of iterations, the graphs in Fig.~\ref{fig:iterations} clearly show that for large numbers of processes $P$, the number of iterations $n$ stays well below the number of processes used, and as a consequence we can expect a speed-up in the (wall) time to solution when using a sufficiently large number of processes.

\subsection{Run times for the algorithms}
\label{sec:evaluation-runtime}

\begin{figure}
    \centering
    \includegraphics[width=0.98\linewidth]{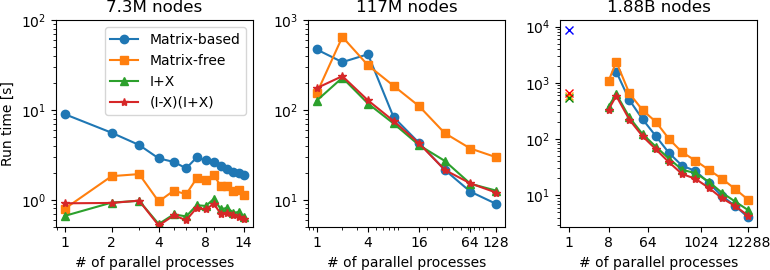}
    \caption{Run times required to solve the flow routing problem. The three panels show results for digital elevation models of the Colorado topography with \num{7.3e6}, \num{117e6}, and \num{1.88e9} node points, respectively solved on a laptop with an Intel Core i9-13900H processor (6 performance plus 8 efficiency cores); a workstation with two AMD EPYC 7773X processors (a total of 128 physical cores); and a cluster with two Intel Xeon Platinum 8160 (``Skylake'') processors per node (for a total of 48 cores per node). Note the different vertical scales.
     For the rightmost panel, computations are too large to be run with $P<12$ (matrix-based) or $P<8$ processes (for the three matrix-free variants); the crosses shown for $P=1$ are estimates obtained by extrapolating from $P=8$ or $P=12$, taking into account how many iterations were needed there and knowing that for $P=1$ one would need only one iteration.}
    \label{fig:runtimes}
\end{figure}

Given the considerations of the previous section, we expect that at least if sufficiently many processes $P$ are used, the solution of the flow routing problem should be substantially faster than on a single process. Fig.~\ref{fig:runtimes} assesses this claim on the three test cases, and using three different machines. Indeed, one can draw a number of conclusions from the data shown there:
\begin{enumerate}
    \item For sufficiently large process counts, say 20 or more, parallel flow routing is indeed \textit{substantially} faster than the computations on a single process. For example, in the middle panel, the fastest method on 128 processes is 14 times faster than the fastest method on a single process (9.0 vs.~127 seconds). The largest model in the right panel can only be solved with $P\ge 8$ processes because of memory restrictions, but comparing the rightmost data point to the estimate provided for $P=1$, it is clear that we would obtain much more significant speed-ups; the speed-up between the largest and smallest problem size computed with the best method (the $I+X$ method) is more than 100; compared to the estimate for $P=1$, the expected speed-up for the matrix-based method would be approximately 2500. In absolute terms, being able to solve a flow routing problem with 1.88 billion points in 4.0 seconds is clearly impressive.
    
    \item With the exception of the matrix-based method, for small numbers of processes, parallel computations cannot compete with the classic high-to-low algorithm of Section~\ref{sec:high-to-low}. This is because for small $P$, the number of Richardson iterations satisfies $n>P$, and so the expected run time $\frac{n}{P}$ (in multiples of the high-to-low algorithm on one process) is greater than one -- i.e., a slow-down. Fig.~\ref{fig:runtimes} indeed bears that out. 
    
    \item For the matrix-based method, the cost of setting up and assembling the matrix (which is independent of the number of Richardson iterations and moreover parallelizes well) substantially exceeds the cost of actually solving the linear system (which does depend on and grows with the number of Richardson iterations). As a consequence, for this method, even small numbers of processes can help reduce the run time.
    
    \item Comparing the three panels of Fig.~\ref{fig:runtimes}, computed on three different machines, it is clear that each machine is better suited to one or the other of the methods. For example, the laptop on which the data for the leftmost panel was generated, clearly does better with the matrix-free methods than the matrix-based one. This is presumably because the processor (an Intel Core i9-13900H) has quite a small cache: 24 MB last-level cache per processor, shared between all of its cores. Consequently, the matrix-based method that requires storing much more data in memory to perform matrix-vector products than the matrix-free ones is severely affected by how slow main memory is compared to cache memory. On the other hand, the workstation used for the middle panel (an AMD EPYC 7773X) has 768 MB for each of the two processors per 128-core machine, and so can store most of the matrix in its caches, substantially accelerating the matrix-based compared to the matrix-free methods. The Xeon Platinum 8160 of the rightmost panel falls in a middle ground with 33 MB of cache for each of the two processors per 48-core machine.
    
    Conversely, the laptop processor has a much higher clock frequency (typically 4.1 GHz, with boost up to 5.4 GHz) compared to the other two machines (typically 2.2 GHz, with boost up to 3.5 GHz for the workstation; and 2.1 GHz with boost up to 3.7 GHz for the cluster nodes) and so, in absolute terms, is substantially faster.

    \item Finally, it may be of interest to compare the results herein with those of \cite{Richardson2014}. To this end, note that Fig.~5 of \cite{Richardson2014} shows results for a case with 182M unknowns and up to 192 processes. Both numbers are approximately 1.5 times those of the middle panel of Fig.~\ref{fig:runtimes}, and so one would expect run times to be about the same. In fact, \cite{Richardson2014} report a run time of approximately a hundred seconds, whereas our program solves the problem in about nine.
\end{enumerate}

The results of this section suggest that the methods we have described herein can equally be applied to much larger problems as well, as long as sufficiently many processes are available.

\section{Conclusions and outlook}
\label{sec:conclusions}

The methods described herein build on previous work by \cite{Richardson2014} but provide substantial improvements using insight about how solvers and preconditioners work and the fact that we have an exact per-process solver available (namely, the high-to-low procedure) that can serve as a very good preconditioner. The results shown in Section~\ref{sec:evaluation} illustrate the performance of these methods on very large and realistic data sets with up to 1.88 billion nodes, and make clear that our algorithms scale well to many thousands of MPI processes.

From the perspective of algorithm design, an algorithm that ``only'' achieves a speed-up of 52 on 128 processes (like the matrix-based method on the medium-sized problem, see the center panel of Fig.~\ref{sec:evaluation-runtime}) is disappointing. Ideally, one would hope to see run times that decrease proportional to $\frac{1}{P}$, but that is not achievable here because the number of iterations \textit{increases} with $P$. At the same time, the following two arguments show that the methods developed herein clearly \textit{are} successes given the original motivation for this work laid out in Section~\ref{sec:introduction}:
\begin{enumerate}
    \item Our goal was to support flow routing on \textit{very large data sets}. These data sets can often not even be loaded on a single process -- as in the case of the large model --, but they do fit into memory when one has hundreds or thousands of machines available. All of the algorithms described herein almost perfectly distribute memory, i.e., if it takes $M$ bytes to do the work on one process, then the methods we describe require $\frac{M}{P}$ bytes on each process to store the local routing data, each process's matrix rows and vector components, etc., perhaps plus a small amount of memory related to the fact that we also need to know information about one layer of nodes around the perimeter of the process's own nodes. $M$ bytes may be too large to store on one machine, but we can often choose $P$ large enough to make $\frac{M}{P}$ bytes fit on every machine of a cluster.

    \item Moreover, the best-performing methods are all ``matrix-free'' in which the action of the operators are described by functions that, as input, only require the local routing information. This information is extremely compact: We only need to store a list of two integers for each node, describing pairs (\textit{src}$\rightarrow$\textit{dst}) of who gives water to whom, sorted high to low with regard to the elevation of the \textit{src} node. Not only is this data structure small and compact, it is also cache-efficient and readily available in every programming language.

    \item A speed-up of 52 translates into a reduction of run time from 472 seconds to 9.0 seconds (again considering the matrix-based algorithm on 117 million nodes). If flow routing is a component in a bigger computation -- say, in time stepping the evolution of a landscape in erosion model simulations -- then a reduction by a factor of 52 is often the difference between a component that is a bottleneck of a simulation, and one that is not.

    \item In absolute terms, the methods described herein can solve a flow routing problem with 1.88 billion points in 4.0 seconds. This is clearly impressive.
\end{enumerate}

As a consequence, we believe that the methods described herein are in fact successful and present a substantial advance over the state of the art.

Looking beyond the results of this paper, it may also be interesting to ask where one might go from here. First, it seems reasonable to believe that a good starting guess would substantially reduce the number of Richardson iterations necessary -- in fact, \cite{Richardson2014} describe exactly this in their Figure 13 where they use the previous time step's solution vector as the starting guess for the current time step in a time-dependent simulation. It should be possible to generalize this principle in a multilevel scheme: Compute the solution on a coarse DEM, and then use this solution as the starting guess on a finer DEM, perhaps recursively. Because solving on a coarser DEM is much cheaper, the effort there will have little impact on the overall cost, but the reduction in iterations on the finer DEM will. The key question one will have to answer is how to create a starting guess on one DEM from the solution on a coarser DEM, a question made more difficult by the fact that the water flow field is not a smooth function of space.

Second, in the computations shown herein, the depression filling step that was necessary before flow routing (see Section~\ref{sec:evaluation}) is a substantial bottleneck. The program that accompanies this manuscript uses the priority-flood algorithm \cite{Barnes2014}, which is sequential. What is needed is a parallel method in which no process has access to the entire DEM. \cite{Barnes2016} presents a parallel depression-filling algorithm; it subdivides the domain into rectangular tiles and scales well to 48 cores. The limitation to tiles (rather than the arbitrarily-shaped subdomains used herein) is likely immaterial, but the use of a single-producer, multiple-consumer design (or join-fork) algorithm will need to be revisited to allow the algorithm to scale to the much larger number of processes we use here.

%
%

\section*{Acknowledgments}

I appreciate the insights of my colleagues on the Landlab+ASPECT team (Mark Behn, Daniel Douglas, Tian Gan, Timo Heister, Eric Hutton, Leif Karlstrom, John Naliboff, Greg Tucker) who first taught me about flow routing, the traditional algorithm described in Section~\ref{sec:high-to-low}, and the difficulty in parallelizing it. They also greatly helped with getting an overview of the existing literature.

I would like to specifically acknowledge an off-hand comment by Timo Heister who suggested that one could perhaps come up with a matrix-free version of the algorithm in Section~\ref{sec:basic-algorithm}, as shown in Section~\ref{sec:optimization-matrix-free}. Over the course of a long, lonely bike ride a few weeks later, my mind morphed his comment into the material of Sections~\ref{sec:optimization-combined-I+X} and \ref{sec:optimization-I+X-times-I-X} -- almost certainly not what he had had in mind at the time, but he gets credit for planting a seed!

The author's work was partially supported by the National Science Foundation under awards EAR-1925595 and OAC-2410847. The computations shown in the right panels of Figures~\ref{fig:iterations} and \ref{fig:runtimes} were performed on the Stampede3 cluster at the University of Texas, Austin, a machine supported by the National Science Foundation. The allocation there is greatly appreciated.

\appendix
\section{Extension to multiple flow direction (MFD) routing schemes}
\label{sec:Dinfty}

The algorithms outlined in the main part of the paper have all assumed that each node gives its water to exactly one of its neighbors -- in the case of the $D_4$ scheme, one among its east, west, north, or south neighbors; in the case of the $D_8$ scheme also perhaps one of the diagonal neighbors; on any kind of mesh, one would typically use the neighboring node to which the gradient from the current node is the most negative.

On the other hand, there are ``multiple flow direction'' (MFD) schemes such as the the ones discussed in \cite{Freeman1984} or \cite{Tarboton1997} in which the water from each node is \textit{partitioned} among its lower-lying neighbors; \cite{Li2022} contains an excellent discussion of the many schemes of this kind that have been proposed. In that case, the algorithms outlined in Section~\ref{sec:better-algorithm} do not immediately apply. However, they are easily generalized.

The starting point of the generalization is that instead of a list of pairs (\textit{src}$\rightarrow$\textit{dst}), the local flow routing information is given by a list of objects of the form 
\begin{align*}
  (\text{\textit{src}}\rightarrow[(dst_1,f_1), (dst_1,f_2), \ldots]),    
\end{align*}
where $dst_n$ is the $n$th node receiving water from $src$ and $f_n$ is the fraction of $src$'s water $dst_n$ receives. As before, one easily convinces oneself that each of these objects corresponds to one column of the matrix, and so the generalization of the matrix-vector product function from Section~\ref{sec:optimization-matrix-free} would now look as follows:
\begin{codebox}[]{matrix-vector product with $A$; multiple water receivers}
\begin{lstlisting}
def apply_A (x):
  y = x                     
  for src,dsts_and_fractions in local_flow_routing_list:
    for dst,f in dsts_and_fractions:
      y(dst) = y(dst) - f*x(src)
  return y
\end{lstlisting}
\end{codebox}
Note that we no longer have to check whether the \textit{dst} element of a pair is invalid, as in Section~\ref{sec:optimization-matrix-free}: If a source node has no destination node, then the list that forms the second part of the \texttt{local\_flow\_routing\_list} objects is simply empty, and the inner loop never executes.

It is clear how the other functions outlined in Section~\ref{sec:optimized-algorithm} generalize in the same way to the case where each node may give water to more than one neighbor.

\section{Implementing the algorithms in practice}
\label{sec:implementation}

An important consideration in \cite{Richardson2014} was that the algorithms described therein can be implemented based on one of the large and widely-used parallel linear algebra packages -- in their case, on the PETSc library \cite{petsc-user-ref,petsc-web-page}. That mattered in order to test the many solver and preconditioner options provided by PETSc. In a similar vein, the program for the current manuscript is implemented using another widely used library, namely \dealii{} \cite{arndt2019dealii,Arndt2025}, using PETSc for the underlying linear algebra.

At the same time, having figured out in Sections~\ref{sec:optimized-algorithm} and \ref{sec:evaluation} that matrix-free algorithms can effectively solve the problem, it turns out that no dedicated linear algebra library is actually necessary. This is due to two reasons: (i) We do not need to build matrix or preconditioner objects as it is enough to provide functions that implement the \textit{actions} of these objects; indeed, these functions all have no more than a handful of lines of code and can easily be implemented in any programming language. (ii) The Richardson iteration we use is quite trivial to implement by hand, without the need to build on a large library that may be difficult to learn. In the following two subsections, let us therefore discuss the specifics of the Richardson solver, and what other software pieces are necessary to implement the methods of this paper.

\subsection{The Richardson solver}
\label{sec:richardson}
Recall from Section~\ref{sec:basic-algorithm} that we want to solve the linear system $A\mathbf w=\mathbf r$ using an iterative solver with the matrix $B=D^{-1}$ as preconditioner. Richardson's iteration repeatedly computes
\begin{align}
\label{eq:richardson-A}
    \mathbf w^{(\ell+1)}
    =
    \mathbf w^{(\ell)} + D^{-1}(\mathbf r - A\mathbf w^{(\ell)}),
\end{align}
which only requires storing the previous and current vectors,  $\mathbf w^{(\ell)}$ and $\mathbf w^{(\ell+1)}$, two temporary vectors, and applying the matrix and precondition operations $A$ and $D^{-1}$ for which Section~\ref{sec:optimized-algorithm} provides implementations. The iteration is stopped whenever the vector in parentheses on the right hand side is sufficiently small, i.e., whenever the residual satisfies $\|\mathbf r - A\mathbf w^{(\ell)}\|<\text{tol}$. It is clear that implementing this method is straightforward, even without relying on large existing libraries.

It is perhaps interesting to see how the iteration would look like for the I+X method of Section~\ref{sec:optimization-combined-I+X}. There, we want to solve $(I+X)\mathbf w=D^{-1}\mathbf r$. If we recall that $X=D^{-1}R$, then the Richardson iteration reads as follows:
\begin{align}
\label{eq:richardson-I+X}
    \mathbf w^{(\ell+1)}
    &=
    \mathbf w^{(\ell)} + (D^{-1}\mathbf r - (I+X)\mathbf w^{(\ell)})
    \notag
    \\
    &=
    D^{-1}(\mathbf r+R\mathbf w^{(\ell)}),
\end{align}
which is even simpler. Recall that the iterations \eqref{eq:richardson-A} and \eqref{eq:richardson-I+X} are identical from a mathematical perspective and will lead to the same iterates $\mathbf w^{(\ell)}$ if performed in exact arithmetic.

The representation \eqref{eq:richardson-I+X} also allows for an interpretation of what the algorithm actually does. To this end, recall that $D^{-1}=(D^{pp\downarrow})^{-1}$ corresponds to the action of transporting water from nodes owned by the current parallel process to all the downhill nodes also owned by the current process. Now, first, $\mathbf r$ corresponds to the water that falls onto the areas that surround each node. Second, $R$ is a matrix that only consists of the \textit{off-diagonal} blocks of $A$; applying $R$ to the vector $\mathbf w^{(\ell)}$ corresponds to the action of passing the water from nodes \textit{not} owned by the current process but adjacent to it, to nodes owned by the current process. We can then think of $\mathbf r+R\mathbf w^{(\ell)}$ as a vector that describes the amount of water each node owned by the current process receives either through rain or from neighboring uphill off-process nodes. The subsequent application of $D^{-1}$ then transports this water downstream.

As a consequence, if the amount of water $\mathbf w^{(\ell)}$ is correct for all the upstream nodes \textit{surrounding} the current process's subdomain, then $\mathbf w^{(\ell+1)}$ will be correct for all the nodes \textit{within} the current process's subdomain. It is easy to see that if the domain contains streams that cross boundaries between subdomains $C$ times, then it takes $C+1$ iterations to get the amount of water correct on all processes that have part of this stream: The first iteration from the headwaters to the first subdomain boundary crossing; the second iteration from that point to the second subdomain crossing; and so on. The argument also supports the claim made at the end of Section~\ref{sec:better-algorithm} that we should expect the number of iterations to increase as we add more and more processes: Because splitting the domain into smaller and smaller chunks must surely increase the number of subdomain crossings $C$ streams have to undergo. On the other hand, it is not \textit{necessary} to run as many iterations as there are subdomain crossings: The solution on a stream with $C$ subdomain crossings will be \textit{exact} after $C+1$ iterations, but in many cases, the residual $\|\mathbf r - A\mathbf w^{(\ell)}\|$ will be small enough for practical purposes to terminate the iteration earlier.

Along similar lines, one can understand the $(I-X)(I+X)$ scheme of Section~\ref{sec:optimization-I+X-times-I-X} as a variation in which one takes the equivalent of \textit{two} of the steps above within each Richardson iteration. It is, then, not a surprise to see that the number of iterations for this scheme is reduced by pretty much exactly a factor of two, see Fig.~\ref{fig:iterations}, compared to the $(I+X)$ scheme of Section~\ref{sec:optimization-combined-I+X}.

\subsection{Software components necessary for parallel flow routing}

As discussed above, the methods described herein are simple enough not to require external parallel linear algebra packages: One can get away with only storing vectors, implementing the application of matrices $A,R,D^{-1}$ in the form of \textit{functions}, and implementing the Richardson iteration itself. As input, these functions only require knowledge of the list of (\textit{src}$\rightarrow$\textit{dst}) pairs that results from the local flow routing step and which is both easy to compute and store. There are, however, two non-trivial pieces for which the use of external tools is useful for parallel implementations: (i) partitioning the domain, and (ii) the ``ghost exchange'' of vector elements; I will detail both briefly in the following.

Conceptually, a digital elevation model (DEM) can be seen as a set of nodes that are connected by edges if we want to consider routing water between two nodes. In other words, a DEM is a \textit{graph}. This is true if the DEM is given as a rectangular grid of nodes (perhaps at fixed latitudes and longitudes, or in an $x$-$y$ coordinate system), but also if it is an unstructured finite element mesh. For the parallel algorithms herein, we need to \textit{partition} this graph into subdomains, preferably without requiring a single node to have to store the entire DEM because we expect no node to have enough memory to store the entire DEM. This can be facilitated in a variety of ways. For the work herein, we use the \dealii{} library \cite{arndt2019dealii,Arndt2025} that internally uses \texttt{p4est} \cite{burstedde2011p4est} for mesh partitioning. Another widely used option is the ParMETIS package \cite{karypis1998metis}. All commonly used partitioning packages then present each parallel process with that sub-graph that is owned by that process, and if necessary one or more layers of nodes owned by other processes that are connected by an edge to the locally owned sub-graph -- so-called ``ghost nodes''. This is enough, then, for each process to request the elevation data from the DEM for all of its own nodes as well as ghost nodes. If the DEM is stored in formats such as NetCDF or HDF5, supporting libraries can then provide this elevation data without having to load \textit{all} of the data.

The second component is that one needs to store the $\mathbf w$ and $\mathbf r$ vectors, along with temporary vectors for the operations in \eqref{eq:richardson-A} or \eqref{eq:richardson-I+X}. One should think of these as \textit{global} vectors, i.e., as a single vector defined on the graph of nodes, where each process then only stores those vector entries that belong to nodes owned by that process plus perhaps the values that correspond to that process's ghost nodes. The vector entries that correspond to ghost nodes are called ``ghost entries''; they may be stored on a process, but are \textit{owned} by another. Most of the operations discussed in previous sections only operate on a process's locally owned vector entries; this is particularly true when multiplying a vector by $D^{-1}$. On the other hand, multiplying by the $R$ matrix (or by $X=D^{-1}R$) requires reading ghost entries. Many libraries -- including \dealii{} and PETSc -- provide data structures for parallel vectors that also allow access to ghost entries and that allow a process to import the values of ghost entries from their owning processes.

%
%

\bibliographystyle{cas-model2-names}
\bibliography{paper}

@article{Richardson2014,
author = {Richardson, Alan and Hill, Christopher N. and Perron, J. Taylor},
title = {{IDA}: An implicit, parallelizable method for calculating drainage area},
journal = {Water Resources Research},
volume = {50},
number = {5},
pages = {4110-4130},
doi = {10.1002/2013WR014326},
url = {https://agupubs.onlinelibrary.wiley.com/doi/abs/10.1002/2013WR014326},
eprint = {https://agupubs.onlinelibrary.wiley.com/doi/pdf/10.1002/2013WR014326},
year = {2014}
}

@article{Arndt2025,
  author = { Arndt, Daniel and Bangerth, Wolfgang and Bergbauer, Maximilian and Blais, Bruno and Fehling, Marc and Gassm\"oller, Rene and Heister, Timo and Heltai, Luca and Kronbichler, Martin and Maier, Matthias and Munch, Peter and Scheuerman, Sam and Turcksin, Bruno and Uzunbajakau, Siarhei and Wells, David and Wichrowski, Michał },
  title = { The {deal.II} library, version 9.7 },
  year = { 2025 },
  journal = {Journal of Numerical Mathematics},
  volume = { 33 },
  issue = { 4 },
  pages = { 403--415 },
  doi = {10.1515/jnma-2025-0115}
}

@Article{arndt2019dealii,
  Title                    = {The {deal.II} finite element library: design, features, and insights},
  Author                   = {D. Arndt and W. Bangerth and D. Davydov and T. Heister and L. Heltai and M. Kronbichler and M. Maier and J.-P. Pelteret and B. Turcksin and D. Wells},
  Journal                  = {Computers \& Mathematics with Applications},
  Year                     = {2021},
  Pages                    = {407-422},
  Volume                   = {81},
  Doi                      = {10.1016/j.camwa.2020.02.022},
  ISSN                     = {0898-1221},
  Url                      = {https://arxiv.org/abs/1910.13247}
}

@Article{BBHK11,
  author =       {W. Bangerth and C. Burstedde and T. Heister
                  and M. Kronbichler},
  title =        {Algorithms and data structures for massively parallel generic
  adaptive finite element codes},
  journal =      {ACM Transactions on Mathematical Software},
  year =         2011,
  volume =       38,
  pages =        {14/1--28}}

@article{Tarboton1997,
  author = { Tarboton, David G. },
  title = { A new method for the determination of flow directions and upslope areas in grid digital elevation models },
  journal = { Water Resources Research },
  year = { 1997 },
  volume = { 33 },
  issue = { 2 },
  pages = { 309--319 },
  doi = {10.1029/96WR03137},
  url = {http://doi.org/10.1029/96WR03137},
}

@article{Freeman1984,
title = {Calculating catchment area with divergent flow based on a regular grid},
journal = {Computers \& Geosciences},
volume = {17},
number = {3},
pages = {413-422},
year = {1991},
issn = {0098-3004},
doi = {https://doi.org/10.1016/0098-3004(91)90048-I},
url = {https://www.sciencedirect.com/science/article/pii/009830049190048I},
author = {T.Graham Freeman}
}

@article{OCallaghanMark1984,
title = {The extraction of drainage networks from digital elevation data},
journal = {Computer Vision, Graphics, and Image Processing},
volume = {28},
number = {3},
pages = {323-344},
year = {1984},
issn = {0734-189X},
doi = {https://doi.org/10.1016/S0734-189X(84)80011-0},
url = {https://www.sciencedirect.com/science/article/pii/S0734189X84800110},
author = {John F. O'Callaghan and David M. Mark}
}

@article{Jenson1988,
  title={Extracting topographic structure from digital elevation data for geographic information system analysis},
  author={Jenson, Susan K and Domingue, Julia O},
  journal={Photogrammetric engineering and remote sensing},
  volume={54},
  number={11},
  pages={1593--1600},
  year={1988},
  publisher={Maryland}
}

@misc{CopernicusData,
  author = {{European Space Agency}},
  title = {{Copernicus Global Digital Elevation Model, distributed by OpenTopography}},
  howpublished = {\url{https://doi.org/10.5069/G9028PQB}},
  year = {2026},
  note = {Accessed: 2026-03-19},
  doi = {10.5069/G9028PQB}
}

@article{Li2022,
author = {Li, Zhenya and Shi, Pengfei and Yang, Tao and Wang, Chao and Yong, Bin and Song, Ying},
title = {An Improved {D8-LTD} for the Extraction of Total Contributing Area ({TCA}) by Adopting the Strategies of Path Independency and Local Dispersion},
journal = {Water Resources Research},
volume = {58},
number = {2},
pages = {e2021WR030948},
doi = {https://doi.org/10.1029/2021WR030948},
url = {https://agupubs.onlinelibrary.wiley.com/doi/abs/10.1029/2021WR030948},
eprint = {https://agupubs.onlinelibrary.wiley.com/doi/pdf/10.1029/2021WR030948},
year = {2022}
}

@article{Barnes2014,
author = {R. Barnes and C. Lehman and D. Mulla},
title = {Priority-flood: An optimal depression-filling and watershed-labeling algorithm for digital elevation models},
journal = {Computers \& Geosciences},
volume = {62},
pages = {117-127},
year = {2014},
issn = {0098-3004},
doi = {https://doi.org/10.1016/j.cageo.2013.04.024},
url = {https://www.sciencedirect.com/science/article/pii/S0098300413001337}
}

@Article{burstedde2011p4est,
  author =       {C. Burstedde and L. C. Wilcox and O. Ghattas},
  title =        {\texttt{p4est}: {S}calable algorithms for parallel
                  adaptive mesh refinement on forests of octrees},
  journal =      {SIAM J. Sci. Comput.},
  volume =       33,
  number =       3,
  pages =        {1103-1133},
  year =         2011}

@article{karypis1998metis,
  title={A fast and high quality multilevel scheme for partitioning irregular graphs},
  author={Karypis, G. and Kumar, V.},
  journal={SIAM J. Sci. Comput.},
  volume={20},
  number={1},
  pages={359--392},
  year={1998},
  publisher={SIAM}
}

@Misc{            petsc-web-page,
  author        = {Satish Balay and Shrirang Abhyankar and Mark~F. Adams and Steven Benson and Jed
                  Brown and Peter Brune and Kris Buschelman and Emil~M. Constantinescu and Lisandro
                  Dalcin and Alp Dener and Victor Eijkhout and Jacob Faibussowitsch and William~D.
                  Gropp and V\'{a}clav Hapla and Tobin Isaac and Pierre Jolivet and Dmitry Karpeev
                  and Dinesh Kaushik and Matthew~G. Knepley and Fande Kong and Scott Kruger and
                  Dave~A. May and Lois Curfman McInnes and Richard Tran Mills and Lawrence Mitchell
                  and Todd Munson and Jose~E. Roman and Karl Rupp and Patrick Sanan and Jason Sarich
                  and Barry~F. Smith and Stefano Zampini and Hong Zhang and Hong Zhang and Junchao
                  Zhang},
  title         = {{PETS}c {W}eb page},
  url           = {https://petsc.org/},
  howpublished  = {\url{https://petsc.org/}},
  year          = {2026}
}

@TechReport{      petsc-user-ref,
  author        = {Satish Balay and Shrirang Abhyankar and Mark~F. Adams and Steven Benson and Jed
                  Brown and Peter Brune and Kris Buschelman and Emil Constantinescu and Lisandro
                  Dalcin and Alp Dener and Victor Eijkhout and Jacob Faibussowitsch and William~D.
                  Gropp and V\'{a}clav Hapla and Tobin Isaac and Pierre Jolivet and Dmitry Karpeev
                  and Dinesh Kaushik and Matthew~G. Knepley and Fande Kong and Scott Kruger and
                  Dave~A. May and Lois Curfman McInnes and Richard Tran Mills and Lawrence Mitchell
                  and Todd Munson and Jose~E. Roman and Karl Rupp and Patrick Sanan and Jason Sarich
                  and Barry~F. Smith and Hansol Suh and Stefano Zampini and Hong Zhang and Hong Zhang
                  and Junchao Zhang},
  title         = {{PETSc/TAO} Users Manual},
  institution   = {Argonne National Laboratory},
  number        = {ANL-21/39 - Revision 3.24},
  doi           = {10.2172/2998643},
  year          = {2025}
}

@misc{ASTER,
  author = {{NASA/METI/AIST/Japan Spacesystems and U.S./Japan ASTER Science Team}},
  title = {{ASTER Global Digital Elevation Model V003}},
  howpublished = {\url{https://www.earthdata.nasa.gov/data/catalog/lpcloud-astgtm-003}},
  year = {2019},
  note = {Accessed: 2026-03-19},
  doi = {10.5067/ASTER/ASTGTM.003}
}

@article{Schwanghart2010,
author = {Wolfgang Schwanghart and Nikolaus J. Kuhn},
title = {TopoToolbox: A set of Matlab functions for topographic analysis},
journal = {Environmental Modelling \& Software},
volume = {25},
number = {6},
pages = {770-781},
year = {2010},
issn = {1364-8152},
doi = {https://doi.org/10.1016/j.envsoft.2009.12.002},
url = {https://www.sciencedirect.com/science/article/pii/S1364815209003053}
}

@misc{Eddins2007,
  author = {Eddins, S.},
  title = {Upslope area—Forming and solving the flow matrix, MathWorks.},
  howpublished = {\url{http://blogs.mathworks.com/steve/2007/
08/07/upslope-area-flow-matrix/}},
  year = {2007}
}

@inproceedings{wallace2010parallel,
  title={Parallel Algorithms for Processing Hydrologic Properties from Digital Terrain},
  author={Wallace, R.M. and Tarboton, D.G. and Watson, D.W. and Schreuders, K.A.T. and Tesfa, T.K.},
  booktitle={Proceedings of the 6th International Conference on Geographic Information Science (GIScience 2010)},
  year={2010},
  pages = {2540--2545},
  address={Zurich, Switzerland},
  url={https://giscience2010.org/pdfs/paper_229.pdf}
}

@article{Braun2013,
title = {A very efficient O(n), implicit and parallel method to solve the stream power equation governing fluvial incision and landscape evolution},
journal = {Geomorphology},
volume = {180-181},
pages = {170-179},
year = {2013},
issn = {0169-555X},
doi = {https://doi.org/10.1016/j.geomorph.2012.10.008},
url = {https://www.sciencedirect.com/science/article/pii/S0169555X12004618},
author = {Jean Braun and Sean D. Willett}
}

@article{Richardson1911,
  author = { Richardson, Lewis Fry },
  title = {The approximate arithmetical solution by finite differences of physical problems involving differential equations, with an application to the stresses in a masonry dam },
  journal = {Philosophical Transactions of the Royal Society of London, Series A: Containing Papers of a Mathematical or Physical Character},
  year = { 1911 },
  volume = { 210 },
  issue = { 459-470 },
  pages = { 307--357 },
  doi = {10.1098/rsta.1911.0009},
  url = {http://doi.org/10.1098/rsta.1911.0009},
  eprint = {https://royalsocietypublishing.org/rsta/article-pdf/210/459-470/307/264773/rsta.1911.0009.pdf}
}

@Book{QV99,
  author =       {A. Quarteroni and A. Valli},
  title =        {Domain Decomposition Methods for Partial Differential Equations},
  publisher =    {Clarendon Press, Oxford},
  year =         1999
}

@inproceedings{Hysom99,
author = {Hysom, David and Pothen, Alex},
title = {Efficient parallel computation of ILU(k) preconditioners},
year = {1999},
isbn = {1581130910},
publisher = {Association for Computing Machinery},
address = {New York, NY, USA},
url = {https://doi.org/10.1145/331532.331561},
doi = {10.1145/331532.331561},
booktitle = {Proceedings of the 1999 ACM/IEEE Conference on Supercomputing},
pages = {29–es},
location = {Portland, Oregon, USA},
series = {SC '99}
}

@article{Hysom01,
author = {Hysom, David and Pothen, Alex},
title = {A Scalable Parallel Algorithm for Incomplete Factor Preconditioning},
journal = {SIAM Journal on Scientific Computing},
volume = {22},
number = {6},
pages = {2194-2215},
year = {2001},
doi = {10.1137/S1064827500376193}
}

@article{Barnes2016,
title = {Parallel Priority-Flood depression filling for trillion cell digital elevation models on desktops or clusters},
journal = {Computers \& Geosciences},
volume = {96},
pages = {56-68},
year = {2016},
issn = {0098-3004},
doi = {https://doi.org/10.1016/j.cageo.2016.07.001},
url = {https://www.sciencedirect.com/science/article/pii/S0098300416301704},
author = {Richard Barnes}
}

@article{Arge2003,
  author = { Arge, Lars and Chase, Jeffrey S. and Halpin, Patrick and Toma, Laura and Vitter, Jeffrey S. and Urban, Dean and Wickremesinghe, Rajiv },
  title = { Efficient Flow Computation on Massive Grid Terrain Datasets },
  journal = { GeoInformatica },
  year = { 2003 },
  volume = { 7 },
  issue = { 4 },
  pages = { 283--313 },
  doi = {10.1023/A:1025526421410},
  url = {http://doi.org/10.1023/A:1025526421410},
}

\end{document}